\documentclass[a4paper,twoside,11pt]{amsart}

\usepackage{amsthm,amsmath,amssymb,amscd}
\usepackage[greek,english]{babel}
\usepackage{graphicx}

\newtheorem{theorem}{Theorem}[section]
\newtheorem{proposition}[theorem]{Proposition}
\newtheorem{lemma}[theorem]{Lemma}
\newtheorem{corollary}[theorem]{Corollary}
\newtheorem{Alg}[theorem]{Algorithm}
\theoremstyle{definition}
\newtheorem{definition}[theorem]{Definition}

\newtheorem{example}[theorem]{Example}
\newtheorem{notation}[theorem]{Notation}

\theoremstyle{remark}
\newtheorem{remark}[theorem]{Remark}

\usepackage{epsfig}

\usepackage[latin1]{inputenc} 
\usepackage[T1]{fontenc}

\def\a{\alpha}

\def\b{\beta}
\def\C{\mathbf C}
\def\D{\Delta}

\def\f{\phi}

\def\g{\gamma}

\def\ord{\mbox{\rm ord}}

\def\p{\pi}

\def\Q{\mathbf Q}
\def\R{\mathbf R}

\def\s{\sigma}
\def\t{\tau}
\def\tr{{\underline{t}}}

\def\w{\omega}

\def\Z{\mathbf Z}

\title{Approximate roots, toric resolutions
 and deformations of a plane branch}
\author{P.D. Gonz\'alez P\'erez}

\address{Instituto de Ciencias Matem\'aticas-CSIC-UAM-UC3M-UCM. Departamento de \'Algebra. Facultad de Ciencias Matem\'aticas. Universidad Complutense de Madrid.
Plaza de las Ciencias 3. 28040. Madrid. Spain.}

\email{pgonzalez@mat.ucm.es}


\thanks{Supported by {\em Programa Ram\'on y Cajal} and by  MTM2007-6798-C02-02 grants of {\em Ministerio de Educaci\'on y Ciencia}, Spain.}

\begin{document}

 \begin{abstract}
We analyze the expansions in terms of the approximate roots of  a
Weierstrass polynomial
 $f \in \mathbf{C} \{ x \} [y]$,  defining a plane branch $(C,0)$,
in the light of the toric embedded resolution of the branch.
                                This leads to the
definition of
a class of
(non-equisingular) deformations of a plane branch  $(C,0)$  supported on certain
 monomials in the approximate roots of $f$, which are essential in the study of Harnack smoothings of
real plane branches by Risler and the author. Our results provide
also  a geometrical approach to
 Abhyankar's irreducibility criterion for power series in two
 variables and also a criterion to determine if a family of plane curves
 is equisingular to a
plane branch.
\end{abstract}

\maketitle

\section*{Introduction}

        The use of {\em approximate roots} in the study of plane algebraic curves,
        initiated by
Abhyankar and Moh  in \cite{Moh}, was essential in the proof of
the  famous {\em embedding line theorem} in \cite{Moh2}. Let
$(C,0) \subset (\C^2, 0)$ be a germ of analytically irreducible
plane curve, a {\em plane branch} in what follows. Certain {\em
approximate roots} of the Weierstrass polynomial defining $(C,0)$
are {\em semi-roots}, i.e., they define {\em curvettes} at certain
exceptional divisors of the minimal embedded resolution. A'Campo
and Oka  describe the embedded resolution of a plane branch by a
sequence of toric modifications using approximate roots in
\cite{AC} and give topological proofs of some of the results of
Abhyankar and Moh. See \cite{A-Kyoto,PP,Ploski,Assi,Pinkham} for
an introduction to the notion of approximate root and its
applications.

 We consider canonical local coordinates at an infinitely near point of
the toric embedded resolution, which is defined by the strict
transform of a suitable approximate root (or more generally a {\em
semi-root}) and the exceptional divisor.  In Section
\ref{mono-res} we introduce an injective correspondence between
monomials in these coordinates and monomials in the approximate
roots (see Proposition \ref{comb2}). From this natural
correspondence we derive two applications.

The first application, given in Section \ref{k-exp}, is based on
the relations of the expansions in terms of semi-roots and
Abhyankar's {\em straight line condition} for the {\em generalized
Newton polygons} associated to a plane branch. These relations are
better understood by passing through the toric embedded resolution
of the branch (see Theorem \ref{straight} and Corollary
\ref{coro}). In particular, we prove that the generalized Newton
polygons arise precisely  from the Newton polygons of the strict
transform of $(C,0)$ at the infinitely near points of the toric
embedded resolution of $(C,0)$ (see Remark  \ref{gen-aby}). We
have revisited Abhyankar's irreducibility criterion for power
series in two variables (see \cite{Abhyankar}). We give a proof of
Abhyankar's criterion by using the toric geometry tools we have
previously introduced. As an application we obtain an algorithmic
procedure to decide if  family of plane curves is equisingular to
a plane branch (see Algorithm \ref{cri-1}). This procedure
generalizes the criterion given by A'Campo and Oka in \cite{AC}.

The second one is the definition of a class of (non equisingular)
multi-parametric deformations $C_{\tr}$ of the plane branch, which
we call multi-semi-quasi-homogeneous (msqh). We explain its basic
properties in Section \ref{msqh}. The terms appearing in this
deformation are
 monomials in the semi-roots of $f$. The
 deformation may be seen naturally as a deformation of
 Teissier's embedding of the plane branch
 $C$ in a higher dimensional affine space (see \cite{T}). If the deformation $C_{\tr}$ is
{\em generic} the Milnor number of $(C,0)$ is related to the sum
of the Milnor numbers of some curves defined from $C_{\tr}$ at the
infinitely near points of the toric resolution of $(C,0)$ (see
Proposition \ref{milnor}). As a consequence we obtain a formula
for the Milnor number, which can be seen as a geometrical
realization of the delta invariant of the singularity in terms of
this class of deformations. In a recent joint work with Risler we
apply this class of deformations in the study of the topological
types of smoothings of real plane branches with the maximal number
of connected components (see \cite{GP-R}).

The paper is organized as follows: Section \ref{branches}
introduce basic results and definitions. Section \ref{msqh} only
depends on Section \ref{branches} and \ref{mono-res}.

\section{Plane branches, semi-roots and toric resolution} \label{branches}

See  \cite{Z,Wall,PP,T,A-Kyoto,Campillo, Casas,Tcurvas}, for
references on singularities of algebraic or analytic curves.

\begin{notation}
The ring of formal (resp. convergent)  power series in $x, y$ is
denoted by $\C [[ x, y  ]]$ (resp. by $\C \{ x, y \}$). The {\em
Newton polygon}  $\mathcal{N} (h)$ of a non zero series $h = \sum_{i, j} \a_{i,j} x^i
y^j \in \C [[ x, y]]$ is the convex hull of the set $
\bigcup_{\a_{i, j} \ne 0} \{ (i,j ) + \R^2_{\geq 0} \}$. If
$\Lambda \subset \R^2$  the {\em symbolic restriction} of $h$ to
$\Lambda$ is the polynomial $\sum_{ (i,j) \in \Lambda \cap \Z^2  }
\a_{i,j} x^i y^j $.

If $(C_i, 0) \subset (\C^2, 0)$, $i =1, 2$ are plane curve germs
defined by $h_i (x, y) =0$, for $h_i \in \C \{ x, y \}$, we denote
by $(C_1, C_2)_0$ or by $(h_1, h_2)_0$ the intersection
multiplicity $\dim_\C \C \{ x, y \} /(h_1, h_2) $.
\end{notation}

\subsection{Expansions and approximate roots}

Abhyankar and Moh have applied and developped the expansions using
approximate roots in the study of algebraic curves (see for
instance \cite{Moh, Abhyankar, A-expansions, Moh2}). See the
surveys  \cite{PP,Pinkham,AC,Ploski} on the applications of the
approximate roots in  the study of plane curves.

Let $A$ be a integral domain. Let  $H \in A [y]$ be a monic
polynomial in $y$ of degree $\deg H >0$. Any polynomial $F \in A
[y] $ has a unique {\em $H$-adic expansion} of the form:
\begin{equation} \label{varphi}
F = a_s + a_{s-1} H + \cdots + a_1 H^{s-1} + a_0 H ^s,
\end{equation}
where $a_i \in A [y]$, $\deg a_i < \deg H$ and $s= [\deg F  / \deg
H ]$. The symbol $[a] $ denotes the integral part of $a \in \R$.
This expansion is obtained by iterated Euclidean division by $H$
(see \cite{Z}).

\begin{proposition} \label{g-expand}      {\rm (see \cite{A-expansions} and
\cite{PP})} Let $n_1, \dots, n_g$ be integers $
>1$. If $F_1,\dots, F_{g+1} \in A [y] $ are polynomials of
degrees $1, n_1, n_1 n_2, \dots, n_1 \cdots n_g$ respectively,
then any  polynomial $F \in A [y] $ has a unique expansion of
 the form:
        \begin{equation} \label{a-expand}
F = \sum_I  \a_I  {F}_1^{i_1} \cdots {F}_{g}^{i_{g}}
{F}_{g+1}^{i_{g+1}}, \mbox{ with } \a_I \in A,
\end{equation} where the
components of the index $I=(i_1, \dots, i_{g+1}) $ verify that $ 0
\leq i_1 < n_1$, $\dots,$ $0\leq i_{g} < n_g$, $ 0 \leq i_{g+1}
\leq [ \deg_y F / \deg_y F_{g+1} ]$. Moreover, the degrees in $y$
of the terms ${F}_1^{i_1} \cdots {F}_{g+1}^{i_{g+1}}$ are all
distinct.
\end{proposition}
{\em Proof.} Consider the  $F_{g+1}$-adic expansion, of the form
(\ref{varphi}),  of the polynomial $F$. Iterate the procedure by
taking recursively $F_j$-adic expansions of the coefficients
obtained for $j = 1,\dots, g$ in decreasing order. The assertion
of the degrees in $y$ is consequence of the following elementary
property of the sequence of integers $(n_1, \dots, n_g)$  (see
\cite{PP}, proof of Corollary 1.5.4). $\Box$
\begin{remark} \label{injective}
Let $n_1, \dots, n_g$ be integers greather than $1$. We set
\[\mathcal{A}_{g+1} := \{ I  = (i_1, \dots, i_{g+1}) \mid 0 \leq
i_1 , < n_1, \dots, 0 \leq i_g, < n_g, 0 \leq i_{g+1} \}.\] The
map $\mathcal{A}_{g+1} \rightarrow \Z$, given by $I \mapsto q_I:=
i_1 + n_1 i_2 + \cdots + n_1 \cdots n_{g} i_{g+1} $, is injective.
\end{remark}

Suppose that the integral domain $A$ contains $\Q$. Denote by
$\mathcal{B}_m \subset A [y] $ the set of monic polynomials of
degree $m>0$ in $y$. Let $F \in A[y]$ be a monic polynomial of
degree $N$ divisible by $m$. Suppose that $N = m k$ for some
integer $k \geq 1$. The {\em Tschirnhausen operator} $\t_F:
\mathcal{B}_m \rightarrow \mathcal{B}_m$ is defined by $\t_F (H) =
H + \frac{a_1}{k}$ where $a_1$ is the coefficient of $H^{k-1}$ in
the $H$-adic expansion (\ref{varphi}) of $F$ (in this case  notice
that   $s =k$ in (\ref{varphi}) since $\deg H = m$). For instance,
if $m =1$, $H = y$ and $y' :
  = y + \frac{a_1}{N}$, then
 the coefficient of
$(y')^{N -1}$ in the $y'$-expansion of $F$ is zero. Setting $y'
=\t_F (y)$ defines a change of coordinates, which is classically
called the {\em Tschirnhausen transformation}.

\begin{definition}
Let $A$ a domain containing $\Q$. Let $F \in A [y]$ a monic
polynomial of degree $N$ and suppose $N=m k$. An approximate root
$G$ of degree $m$ of the polynomial $F$  is a monic polynomial in
$A [y]$ such that $\deg ( F- G^{k}) < N - m$.
\end{definition}
The approximate root $G$ of degree $m$ of $F$ exists and is
unique. It is determined algorithmically in terms of Euclidean
division of polynomials by: $ G= \t_F \circ \stackrel{(m)}{\cdots}
\circ \t_F (H), \, \forall \,  H \in \mathcal{B}_m$.

\subsection{ Local toric embedded resolution of a plane
branch}\label{toric-res}

In this paper $(C,0)$ denotes a germ of analytically irreducible
plane  curve, a {\em plane branch} for short, defined by an
irreducible element   in the ring $ \C \{ x, y \}$ of germs of
holomorphic functions at the origin of $\C^2$. We recall the
construction of a {\em local toric embedded resolution of
singularities} of the plane branch  $(C,0)$ by a sequence of
monomial maps. For a complete description see \cite{AC}. See
\cite{Oka1,Oka2,LO,Rebeca} for more on toric geometry and plane
curve singularities.

 We define
a sequence of birational monomial maps $\pi_j: Z_{j+1} \rightarrow
Z_j$, where $Z_{j+1}$ is an affine plane $\C^2$ for $j=1, \dots,
g$, such that the composition $\Pi:= \pi_1 \circ \dots \circ
\pi_g$ is a local {\em embedded resolution} of the plane branch
$(C,0)$, that is, $\Pi$ is an isomorphism over $\C^2 \setminus \{
(0,0) \}$ and the {\em strict transform} $C'$ of the plane branch
$C$ (defined as the closure of the pre-image by $\Pi^{-1}$ of the
punctured curve $C \setminus \{ 0 \}$) is a smooth curve on
$Z_{g+1}$ which intersects the exceptional fiber $\Pi_1^{-1} (0)$
 transversally. Notice that the map $\Pi$ is not proper. The map $\Pi$ can be
 seen as an affine chart of certain sequence of blow-ups of points.

We consider local coordinates $(x, y)$ for $(\C^2, 0)$.
 We
say that $y' \in \C \{ x, y \} $ is  {\em good} with respect to
$(C,0)$ and $\{x= 0 \}$ if setting $(x_1, y_1) := (x, y')$ defines
a pair of local coordinates at the origin and the germ $(C,0)$ is
defined by an equation $f=0$ where,
\begin{equation} \label{normal1}
f = ( y_1^{n_1} - \theta_1 x_1^{m_1})^{e_1} + \cdots,
\end{equation}
in such a way that  $\theta_1 \in \C^*$, $\mbox{\rm gcd } (n_1,
m_1)=1$ and the terms which are not written have exponents $(i,
j)$ such that $i n_1 + j m_1
> n_1 m_1 e_1$, i.e., they lie above the compact edge
$ \Gamma_1 := [(0, n_1 e_1), (m_1 e_1, 0)]$ of the Newton polygon
of $f$. Notice that $e_0 := e_1 n_1$ is the intersection
multiplicity of $(C,0)$ with the line $\{ x_1=0 \}$.

 Such a choice
of $y_1$ is not unique. The choice $y_1 := y + \t_f (y)$, defined
by the {\em Tschirnhausen transformation}, is good with respect to
$\{x_1 = 0 \}$ and $(C,0)$. We assume without loss of generality
that $f$ is a Weierstrass polynomial in $y_1$.

The vector $\vec{p}_1 = (n_1, m_1)$ is  orthogonal to $\Gamma_1$
and defines a subdivision of the  positive quadrant $\R^2_{\geq
0}$, which is obtained by adding the ray $\vec{ p}_1 \R_{\geq 0}$.
The quadrant $\R^2_{\geq 0}$ is subdivided in two cones, $\t_i
:=\vec{ e}_i \R_{\geq 0} + \vec{ p}_1 \R_{\geq 0}$ for $i=1, 2$
where $\{ \vec{ e}_1, \vec{ e}_2 \}$ is the canonical basis of
$\Z^2$. We define the {\em minimal regular subdivision} $\Sigma_1$
of $\R^2_{\geq 0}$ which contains the ray $\vec{ p}_1 \R_{\geq 0}$
by adding the rays defined by those integral vectors in $\R^2_{>0
}$, which belong to the boundary of the convex hull of the sets
$(\t_i \cap \Z^2) \backslash \{ 0\}$, for $i=1,2$.
There is a unique cone $\s_1 = \vec{p}_1 \R_{\geq 0} + \vec{q}_1
\R_{\geq 0}$ in the subdivision $\Sigma_1$ such that $\vec{q}_1 =
(c_1, d_1)$ satisfies that:
\begin{equation} \label{ave-i}
c_1 m_1 - d_1 n_1 = 1.
\end{equation}

By convenience we denote $\C^2$ by $Z_1$,  the coordinates $(x,
y)$ by $(x_1, y_1)$ and the origin $0 \in \C^2 = Z_1$ by  $o_1$.
We also denote $f$ by $f^{(1)}$ and $C$ by $C^{(1)}$.
The map $\p_1: Z_2 \rightarrow Z_1$ is defined by
\begin{equation}
\label{mm-i}
\begin{array}{lcl}
x_1 & = & u_2 ^{c_1} x_2^{n_1},
\\
y_1 & =  & u_2^{d_1} x_2^{m_1},
\end{array}
\end{equation}
where $u_2, x_2$ are coordinates in the affine plane $Z_2 :=
\C^2$. The components of the  exceptional fiber $\p_1^{-1} (0)$
are $\{ x_2 = 0 \}$  and $\{ u_2 =0 \}$.
The  pull-back of $C^{(1)}$ by $\pi_1$ is defined by $f \circ \p_1
=0$. The term $f \circ \p_1$ decomposes as:
\begin{equation} \label{strict-j}
f^{(1)} \circ \p_1 =   \mbox{\rm Exc} (f^{(1)}, \p_1) \,
\bar{f}^{(2)} (x_2, u_2), \mbox{ where } \bar{f}^{(2)} (0,0) \ne
0,
\end{equation}
and $\mbox{\rm Exc} (f^{(1)}, \p_1) := y_1^{e_0} \circ \p_1 = u_2
^{ d_1 e_0} x_2^{m_1 e_0}$. The polynomial $\bar{f}^{(2)} (x_2,
u_2)$ (resp. $ \mbox{\rm Exc} (f^{(1)}, \p_1) $)  defines the {\em
strict transform} $C^{(2)}$ of $C^{(1)}$  (resp.  the {\em
exceptional divisor}).
By  Formula (\ref{normal1}) we find that $ \bar{f}^{(2)} (x_2, 0)
=1$, hence the exceptional line $\{u_2 =0 \}$ does not meet the
strict transform.
Since
\[
 \bar{f}^{(2)}  ( 0, u_2) = ( 1  - \theta_1 u_2^{c_1 m_1 - d_1 n_1}) ^{e_1}
\stackrel{\mbox{ (\ref{ave-i})}}{=}  (1- \theta_1 u_2) ^{e_1},
\]
it follows that $\{ x_2 = 0 \}$ is the only component of the
exceptional fiber of $\pi_1$ which intersects the strict transform
$C^{(2)}$ of $C^{(1)}$, precisely at the point $o_2$ with
coordinates $x_2 =0$ and $u_2=\theta_1^{-1}$ and with intersection
multiplicity equal to $e_1$. If $e_1 =1$ then the map $\pi_1$ is a
local embedded resolution of the germ $(C,0)$. If $e_1 > 1$ we
consider a pair of coordinates $(x_2, y_2)$ at the point $o_2$,
with $y_2$ good for $\{x_2 = 0 \}$ and $(C^{(2)}, o_2)$. It
follows that  $C^{(2)}$ is defined by a term, which we call {\em
the strict transform function}, of the form:
\begin{equation}\label{fg}
 f^{(2)} (x_2, y_2) = (y_2^{n_2} - \theta_2
x_2^{m_2} ) ^{e_2} +  \cdots,
\end{equation}
where  $\theta_2 \in \C^*$, $\mbox{\rm gcd} (n_2, m_2)=1$ and the
terms which are not written have exponents $(i, j)$ such that $i
n_2 + j m_2
> n_2 m_2 e_2$. Notice that $e_1 =
e_2 n_2$.

We iterate this procedure defining for $j>2$ a sequence of
monomial birational maps $\pi_{j-1} :Z_{j} \rightarrow Z_{j-1}$,
which are described by replacing the index $1$ by $j-1$ and the
index $2$ by $j$ above. In particular when we refer to a Formula,
like (\ref{ave-i}) at level $j$, we mean after making this
replacement. We denote by $\mbox{\rm Exc} (f^{(1)} , \p_1 \circ
\cdots \circ \p_{j}
 ) $ the {\em exceptional function} defining the exceptional
 divisor of the pull-back of $C$ by $\pi_1 \circ \cdots \circ
 \pi_{j}$.
 Notice that
 \begin{equation} \label{EJ}
 \mbox{\rm Exc} (f^{(1)}, \p_1 \circ \cdots \circ \p_{j} ) =  (y_1^{e_0}
 \circ \p_{1} \circ \cdots \circ \p_{j}) \cdots (y_{j}^{e_{j-1}}
 \circ \p_{j} ).
 \end{equation}
 Since by
construction we have that $e_j | e_{j-1} | \cdots | e_1 |e_0$ (for
$ | $ denoting divides), at some step we reach a first integer $g$
such that $e_g =1$ and then the process stops. The composition
$\pi_1 \circ \dots \circ \pi_g$ is a local  toric embedded
resolution of the germ $(C,0)$.

\begin{remark} \label{Oka}
Given $e_0 = (x_1, f)_0$,
% If $x_1=0$ is not tangent to $C$ then
the sequence of pairs $\{( m_j, n_j)\}_{j=1}^g$   determines and
it is determined by the {\em characteristic pairs}  or the {\em
Puiseux exponents} of the plane branch $(C,0)$, which are obtained
when the line $\{ x_1 =0 \} $ is not tangent to $C$ at the origin
(see \cite{AC} and \cite{Oka1}). These pairs classify the {\em
embedded topological type} of the germ $(C,0) \subset (\C^2, 0)$,
or equivalently its complex {\em equisingularity type}.
\end{remark}

\begin{notation}
 We set $n_0:=1$. We denote by $f'_j$ the approximate root of the polynomial  $f \in \C \{ x_1 \} [y_1]$,
 of degree $n_0 \cdots n_{j-1}$ in $y_1$,  for $j=1, \dots , g$. The integers $n_i$ are those of
Remark \ref{Oka}. We consider the sequence of intersection
multiplicities given by:
 \begin{equation} \label{m-gen2}
 \bar{b}_0 := e_0 = (x, f)_0, \, \bar{b}_j :=  ( f_j', f)_0, \mbox {
 for } j=1, \dots , g.
 \end{equation}
\end{notation}

\begin{definition} \label{semi-root}
A $j^{th}$-semi-root $(C_j,0)$ of $(C,0)$ with respect to  the
line $\{ x_1 =0 \}$, is a germ $(C_j,0)$ of curve such that $(C_j,
C)_0 = \bar{b}_{j}$  and $ ( C_j, x_1)_0 = n_0 \cdots n_{j-1} $,
for $ 0 \leq j \leq g$. We convey that  $C_{g+1} := C $.
The sequence $\{ (C_j, 0) \}_{j=1}^{g+1}$ is called the
characteristic sequence of semi-roots of $(C, 0)$ with respect to
$\{ x_1 =0 \}$.
\end{definition}

        \begin{remark} For simplicity we have defined semi-roots in terms of approximate roots,
i.e., without passing by  Abhyankar and Moh Theorem (\cite{Moh}). For a definition of semi-roots in terms of
Puiseux exponents and related results see \cite{PP}, for instance.
        \end{remark}

\begin{notation}     \label{not-semi}
Let us fix  a sequence of semi-roots $(C_j,0)$ of the plane branch
$(C,0)$ with respect to $\{x_1 =0 \}$,  for $j=1,\dots, g+1$. Each
curve $C_j$ is defined by a Weierstrass polynomial $f_j \in \C\{
x_1 \} [y_1]$ of degree $n_0 \cdots n_{j-1}$, which we call also
{\em semi-root} by a slight abuse of terminology. We will assume
that $f_1 = y_1$ and $f_{g+1} = f$.
\end{notation}

\begin{definition}
Let us fix  $2 \leq j \leq g$. A germ $(D,0) \subset (\C^2,0)$  is
called a {\em $j^{th}$-curvette} for $(C,0)$  and $\{ x_1 = 0 \}$
if it is analytically irreducible and the strict transform of $D$
by $\p_1 \circ \dots \circ \p_{j-1}$ is smooth and intersects
transversally the exceptional divisor $\{ x_j = 0 \}$ at the point
$o_j \in \{ x_j = 0 \}$. The branch $(D,0)$ is a {\em
$j^{th}$-curvette with maximal contact} if in addition the strict
transform of $(D,0)$ by $\p_1 \circ \dots \circ \p_{j-1}$ is
defined by $y_j'=0$ where $y_j'$ is good with respect to $\{ x_j =
0 \}$ and $(C^{(j)},0)$.
\end{definition}

\begin{proposition} \label{mod} {\em (see \cite{Z-inversion}, \cite{AC}, \cite{PP} and
\cite{GP} \S 3.4)}
\begin{enumerate}
\item[(i)] If $C_j$ is a $j^{th}$-semi-root of $(C,0)$ with
respect to $\{ x_1 =0 \}$ then $(C_j,0)$ is a $j^{th}$-curvette
with maximal contact, for $j=2, \dots, g$.
\item[(ii)] We denote by   $C_{2}^{(2)} , \dots, C_g^{(2)},
C_{g+1}^{(2)} = C^{(2)}  $
 the strict transforms  by the monomial map $\p_1$ of the
semi-roots $C_2, \dots, C_g, C_{g+1} = C$ of the plane branch
$(C,0)$. The sequence $C_{2}^{(2)} , \dots,  C_{g+1}^{(2)}$ is a
characteristic sequence of semi-roots of the branch
$(C^{(2)},o_2)$ with with respect to the line $\{ x_2 =0 \}$.
\end{enumerate}
\end{proposition}

\begin{remark}
We will assume in the rest of the paper that the local coordinate
$y_j$, in the local embedded resolution of $(C, 0)$ introduced
above, is  the strict transform function of the semi-root $f_j$,
for $j= 2, \dots,g$ (we can do this by Proposition \ref{mod}).
This implies that $y_j$ is of the form:
\begin{equation} \label{xi-rel}
y_j = 1- \theta_j u_j + x_j R_j (x_j, u_j)  \mbox{ for some }
 R_j \in \C \{ x_j, u_j \}.
\end{equation}
\end{remark}
As a consequence of Proposition \ref{mod} we have the following:
\begin{remark} \label{one-edge} $\,$
\begin{enumerate}
\item[(i)] If $2 \leq j \leq g$ the Newton polygons of $f (x_1,
y_1)$ and of $f_j^{e_{j-1}} (x_1, y_1) $ have only one compact
edge
 $\Gamma_1$, defined in Section \ref{toric-res}, and the symbolic
 restrictions of $f$ and of $f_j^{e_{j-1}}$ coincide on this edge.
\item[(ii)] If $2 <  j \leq g$ similar statement holds for
$f^{(2)} (x_2, y_2)$ and of $(f_j ^{(2)})^{e_{j-1}} (x_2, y_2) $
and $\Gamma_2$.
 \end{enumerate}
\end{remark}

\begin{definition}
The {\em semigroup of the plane branch} $(C, 0)$ is
$ \Lambda_C: =  \{ (f, h)_0 \mid h \in \C \{ x, y \} - (f)
\}$.
\end{definition}
The semigroup $\Lambda_C$ is generated by     the elements in the sequence
(\ref{m-gen2}).
The sequence  (\ref{m-gen2}) is called the {\em characteristic
sequence of generators} of the semigroup $\Lambda_C$ with respect
to the line $\{ x_1=0 \}$. If the line $\{x_1=0\}$ is not tangent
to $C$ at the origin then the set (\ref{m-gen2}) is a minimal set
of generators of the semigroup $\Lambda_C$ and the notation,
$\bar{\b}_j$ instead of $\bar{b}_j$, is the usual one in the
litterature. The semigroup $\Lambda_C$ has the following
properties (see \cite{T}, for instance).
\begin{lemma} \label{tee}
Any $\bar{b} \in \Lambda_C$ has a unique expansion of the form:
\begin{equation} \label{expand-semi}
\bar{b} = \eta_0 \bar{b}_0 + \eta_1 \bar{b}_1 + \cdots + \eta_g \bar{b}_g,
\end{equation}
where $0 \leq \eta_0$ and $0 \leq \eta_j  < n_{j}$,  for
$j=1,\dots, g$. The image of $\bar{b}_j$ in the group $\Z /
(\sum_{i=0}^{j-1} \Z \bar{b}_i ) $ is of order $n_j$. We have
that:
\begin{equation} \label{plane-semi2}
n_j \bar{b}_j \in \Z_{\geq 0} \bar{b}_0 + \cdots + \Z_{\geq 0} \bar{b}_{j-1}
\mbox{ and }  n_j \bar{b}_j <\bar{b}_{j+1}, \mbox{ for } j=1, \dots, g.
\end{equation}
\end{lemma}

The following Proposition states some numerical relations between the sequences
$\{ (n_j, m_j) \} _{j=1}^g$ and
$(\bar{b}_j)_{j=0}^g$ (see for instance   \cite{GP} \S 3.4).
\begin{proposition} \label{mod2}   We have that
\[(x_j, f^{(j)}
)_{o_j} = e_{j-1} = n_j e_j  \mbox{   and  }  (y_j, f^{(j)} )_{o_j} = \bar{b}_j -
n_{j-1} \bar{b}_{j-1} = m_j e_j   \mbox{ for } 1 \leq j \leq g.\]
\end{proposition}

The following proposition shows the relations between the
characteristic sequences of generators of the semigroups of the
plane branch $(C,0)$ and of its semi-root $C_{j+1}$.

\begin{proposition}      \label{A-M}
Let  $C_{j+1}$ be a $(j+1)^{th}$-semiroot of the plane branch
$(C,0)$, for some $j=1, \dots,g$ (see Definition \ref{semi-root}).
The characteristic sequence of the semigroup of the plane branch
$C_{j+1}$ with respect to the line $\{ x_1 =0 \}$ is equal to $
\frac{1}{e_j} \bar{b}_0, \dots, \frac{1}{e_j} \bar{b}_j$, for
$j=1, \dots,g$ (see (\ref{m-gen2})).
\end{proposition}

The {\em normalization map}  $(\C, 0) \rightarrow (C, 0) $ of the
branch $(C,0)$, which is  of the form $\t \mapsto (x_1( \t), y_1(
\t))$,
 may be defined explicitely in terms of a Newton Puiseux
 parametrization of the branch.
If $ h(x_1,y_1) \in \C \{ x_1 \}[y_1] $ defines a plane curve
germ, we have that $(f, h)_0 = \ord_\t ( h (x_1 ( \t) , y_1 ( \t))
$, where $\ord_\t$ denotes the $\t$-adic valuation of the field
$\C ((\t))$ of Laurent series. We abuse the notation by denoting
with the same letter the functions $u_j, x_j $ and $y_j$ and their
images $u_j(\t), x_j (\t) $ and $y_j (\t)$, induced by the
normalization map, in the field  $\C ((\t))$.

\begin{lemma} \label{j-semi}
We have that $\ord_\t( u_j ) =    0$ and $ \ord_\t(\mbox{\rm Exc}
(f, \p_{1} \circ \cdots \circ \p_{j}))  = e_{j-1} \bar{b}_j$  for
$1 \leq j \leq g$.
\end{lemma}
{\em Proof.} Notice that  $\ord_\t(x_1) = (x_1, f)_0= e_1 n_1$ and
$\ord_\t(y_1) = (y_1, f)_0 =e_1 m_1$. We deduce from (\ref{mm-i})
that $u_2 = x_1^{m_1} y_1^{-n_1}$. It follows that $\ord_\t(u_1) =
0$.
The equality $\ord_\t(\mbox{\rm Exc} (f, \p_{1}) )= e_{0}
\bar{b}_1$, follows from formula (\ref{EJ}). We conclude the proof
by an easy induction on $j$ using Proposition \ref{mod2} and
formula (\ref{EJ}). $\Box$

\begin{example} \label{10}
A local embedded resolution of the real plane branch singularity
$(C, 0)$ defined by $F= (y_1^2 -x_1^3 )^3 - x_1^{10} =0$ is as
follows.
The morphism $\p_1$ of the toric resolution is defined by
\[
\begin{array}{lcl}
x_1 & = & u_2 ^{1} x_2^{2},
\\
y_1   & =  & u_2^{1} x_2^{3}.
\end{array}
\]
We have that  $f_2 := y_1^2 - x_1^3 $ is a $2^{nd}$-curvette for
$(C, 0)$ and $ \{ x_1 = 0  \}$. We have $f_2 \circ \p_1 = u_2 ^2 x_2^ 6 ( 1 -
u_2) =u_2 ^2 x_2 ^ 6 y_2$, where $y_2 := 1 -u_2$ defines the
strict transform function of $f_2$, and together with $x_2$
defines local coordinates at the point of intersection $o_2$ with
the exceptional divisor $\{ x_2=0 \}$. Notice in this case that the term
$R_2$ in (\ref{xi-rel}) is zero.
 For $F$ we find that:
\[
\begin{array}{lcl}
F \circ \p_1 &  = & u_2 ^6 x_2 ^{18 } \left( (1 - u_2)^3 - u_2^{4}
x_2 ^{2} \right).
\end{array}
\]
Hence  $\mbox{\rm Exc} (F, \p_1) :=  y_1^6 \circ \p_1 = u_2 ^6 x_2
^{18 }$ is the exceptional function associated to $F$, and $ F
^{(2) } = y_2^3 - (1- y_2)^4 x_2^2 $ is the strict transform
function. Comparing to (\ref{fg}) we see that $e_2= 1$, $n_2 = 3$,
$m_2 =2$ and the restriction to $F ^{(2) } (x_2, y_2)$ to the
compact edge of its local Newton polygon is equal to $y_2^3
-x_2^2$. The map $\p_2: Z_3 \rightarrow Z_2$ is defined by $x_2=
u_3^2  x_3^{3}$ and $x_3 = u_3 x_3^2$. The composition $\p_1 \circ
\p_2$ defines a local embedded resolution of $(C,0)$.
   \end{example}

\section{Monomials in the semi-roots from the embedded resolution}
\label{mono-res}

We keep notations of the previous Section (cf. Notation
\ref{not-semi} and Remark \ref{Oka}). For $2 \leq j \leq g$ we
consider a sequence of integers of the form
\[
0 \leq i_0,  \quad 0 \leq i_1 < n_1, \quad \dots, \quad    0 \leq
i_{j-1} < n_{j-1}, \quad  0 \leq i_{j} < e_{j-1}.
\]
Notice that by
Proposition \ref{g-expand} the term
\begin{equation} \label{eme}
\mathcal{M}= x^{i_0} f_1^{i_1} \, f_2^{i_2} \cdots f_j^{i_j}
\end{equation}
may appear in the
 $ (f_1, f_2, \dots, f_j)$-expansion of $f$.
For any integer $2 \leq j \leq g$ we define below a map which
associates to a monomial of the form,
\[
x_j^r y_j^s, \mbox{ with } 0 \leq r, s < e_{j-1}
\]
a monomial in $x, f_1, \dots, f_j$ of the form (\ref{eme}). We
study conditions for a term of the form (\ref{eme}) to appear in
the $(f_1, \dots,f_j)$ expansion of $f$. We use these ideas to
analyze equisingular (and non equisingular) classes of
deformations of the branch $(C, 0)$ in the following sections.

\begin{remark} \label{units}
To avoid cumbersome notations if $2 \leq j \leq g+1$  we denote
simply by $u_i$ the term $u_i \circ \p_{i} \circ \cdots \circ
\p_{j-1}$, whenever $i < j$ and the integer $j$ is clear from the
context. The function $u_i \circ \p_{i} \circ \cdots \circ
\p_{j-1}$ has an expansion as a series in $\C \{ x_j, y_j \}$ with
non-zero constant term (see  (\ref{xi-rel}) at level $i < j$).
\end{remark}

The following Lemma is an elementary observation which is useful
to motivate our results:
\begin{lemma} \label{lift} Given a monomial $\mathcal{M} = x_1^{i_0} \,
f_1^{i_1} f_2^{i_2} \cdots f_j^{i_j}$ of the form (\ref{eme})
there exists unique integers    $r,s = i_j$
and   $k_2, \dots, k_j$ such that
\begin{equation} \label{ler}
u_2^{k_2} \cdots u_j^{k_j} \, x_j^r  \, y_j^s = ( \mathcal{M}
\circ \pi_1 \circ \cdots \circ \pi_{j-1} ) \,   ( \mbox{\rm Exc}
(f, \p_{1} \circ \dots \circ \p_{j-1} ))^{-1}.
\end{equation}
The integer $r$ depends only on $\mathcal{M}$ and the sequences
$\{ (n_i, m_i) \}_{i=1}^{j-1}$ and $\{ e_i \}_{i=0}^{j-1}$. The
term $u_j^{k_j}  \cdots u_2^{k_2}$ is a unit in $\C \{ x_j,
y_j\}$.
\end{lemma}
{\em Proof.} By formula (\ref{EJ}) and (\ref{mm-i}) we have that
    $ (\mbox{\rm Exc} (f, \p_{1}))^{-1}  (\mathcal{M} \circ \p_1) =
u_2^{k_2} x_2^{i_0'} y_2^{i_2} (f_3^{(2)})^{i_3} \cdots
(f_j^{(2)})^{i_j} $ for some integer $k_2$ where $i_0' = n_1 i_0 +
m_1 ( -e_0 + i_1 + n_1 i_2 + \cdots + n_1 \cdots n_{j-1} i_j)$. By
Remark \ref{units} the term $u_2^{k_2}$ is a unit in the ring $\C
\{ x_2, y_2 \}$. The result is proved if $j=2$. If $j
> 2$ we find that $ (\mbox{\rm Exc} (f, \p_{1} \circ \p_2 ))^{-1}
(\mathcal{M} \circ \p_1 \circ \p_2) =
      u_2^{k_2}  u_3^{k_3}  x_3^{i_0''} y_3^{i_3} (f_4^{(3)})^{i_4} \cdots (f_j^{(3)})^{i_j} $
for some integer $k_3$ where  $i_0''  = n_2 i_0' + m_2 ( -e_1 + i_2 + n_2 i_3 + \cdots + n_2 \cdots n_{j-1} i_j)$.
The assertion follows by an easy induction on $j$.    $\Box$

\begin{remark} Notice that the condition $ r \geq 0 $ is not
guaranted by Lemma \ref{lift}. See Example \ref{11}.
\end{remark}

The following key Proposition shows that given $(r, s) \in
\Z_{\geq 0}$ with $s < e_{j-1}$ there is a unique way to determine
a suitable monomial ${\mathcal{M}}_j(r,s)$ in $x_1$ and the
semi-roots $y_1 =f_1, f_2, \dots, f_j$, such that  the composite
${\mathcal{M}}_j(r,s) \circ \p_1 \circ \dots \circ \p_{j-1}$ is
equal to the product of  the exceptional divisor function
$\mbox{\rm Exc} (f, \p_{1} \circ \dots \circ \p_{j-1} )$ by the
monomial $x_j^r y_j^s$ times a unit in the ring $\C \{ x_j, y_j
\}$.

\begin{proposition} \label{comb2}
Let us fix a real plane branch $(C,0)$ together with a local toric
embedded resolution $\p_1 \circ \dots \circ  \p_{g}$
(cf.~notations of {\rm Section \ref{toric-res}}). If $2 \leq j
\leq g$ and $(r, s) \in \Z^2_{\geq 0}$ with $s < e_{j-1}$ then
there exists unique integers
\begin{equation} \label{lerdos} 0 < i_0, \quad 0 \leq i_1 < n_1, \quad \dots, \quad 0 \leq
i_{j-1} < n_{j-1}, \quad i_j = s, \end{equation} and $k_2, \dots
k_j
> 0$ such that (\ref{ler}) holds.
\end{proposition}
Recall that the integers $c_1, d_1$ are defined by (\ref{ave-i})
in terms of the pair $(m_1, n_1)$.
\begin{lemma} \label{dos}
If $r \geq 0, l > 0$ are integers there exist unique integers $k,
i_0, i_1$  such that
$
u_2^{k} x_2^{r}    = ((x_1^{{i_0}}
y_1^{{i_1}}) \circ \p_1 ) (y_1^{l n_1} \circ \p_1)^{-1} $   with $0< {i_0}, k$
and $ 0
\leq {i_1} < n_1$.
 We have that:
\begin{equation} \label{m-k3}
k = l + [ {c_1} r /{n_1}  ], \quad {i_0} = k m_1 - r d_1  \mbox{
and } \quad {i_1} = c_1 r - n_1 [ {c_1} r /{n_1} ].
\end{equation}
In particular, ${i_1} =0$ if and only if $r= p n_1$ for some
integer $p$.
%  and then $k= l + p c_1 $ and ${i_0} = l m_1 + p$.
\end{lemma}
{\em Proof.} By (\ref{mm-i}) we deduce that $u_2 = x_1^{m_1} y_1
^{-n_1}$ and $x_2 = x_1^{-d_1} y_1 ^{c_1}$. The term \[
 (y_1^{l n_1}\circ \p_1 ) u_2^{k}  x_2 ^{r}   = (x_1^{ k m_1 - r d_1}
y_1^{ r c_1 +  (l-k) n_1}) \circ \p_1 \] is the transform of a
holomorphic monomial by $\p_1$ if and only if: \[  0 \leq {i_0'}
:= { k m_1 - r d_1} \mbox{ and } 0 \leq  {i_1'} := r c_1 + ( l -
k) n_1,\]
 or equivalently, $
\frac{d_1}{m_1} r \leq k \leq \frac{c_1}{n_1} r + l$. By
(\ref{ave-i}) we have that $m_1 c_1 - d_1 n_1 = 1$. This implies
that $\frac{d_1}{m_1} < \frac{c_1}{n_1}$, thus the interval of the
real line $[ \frac{d_1}{m_1} r ,  \frac{c_1}{n_1} r + l ] $ is of
length greater than $l \geq 1$. Any integer $k$ lying on this
interval is convenient to define a holomorphic monomial. The
condition ${i_1'} < n_1$, is equivalent to $\frac{c_1}{n_1} r + l
-k < 1 $, and it is verified if and only if $k = [ \frac{c_1}{n_1}
r + l] = l +  [ \frac{c_1}{n_1} r] >0$. We denote the integers
${i_0'}$ and ${i_1'}$ corresponding to this choice of $k$ by
${i_0}$ and $i_1$ respectively. We have that:
\[
{i_0}  = ( \frac{c_1}{n_1} r + l ) m_1 - r d_1  > (
\frac{c_1}{n_1} r + l -1 ) m_1 - r d_1 = r m_1 (  \frac{c_1}{n_1}
- \frac{d_1}{m_1} ) + (l-1 ) m_1
\stackrel{\mbox{\rm(\ref{ave-i})}}{\geq} ( l -1 ) m_1 \geq 0.
\]

For the last assertion, we have that ${i_1} = c_1 r - n_1 [ {c_1}
r /{n_1} ] =0$ if and only if  $n_1$ divides $r$, since
$\mbox{{\rm gcd}} (c_1, n_1) =1$ by (\ref{mm-i}). $ \Box$

\begin{lemma} \label{k1}
If $(r, s) \in \Z_{\geq 0}$ with $s < e_1$ there exist unique
integers $k, {i_0}, {i_1}$ with $0 < k , {i_0}$ and $0 \leq {i_1}
< n_1$ such that: $  u_2^{k} {x_2}^r
{y_2} ^s = (({x_1}^{i_0} y^{i_1} f_2^s ) \circ \p_1)  ( \mbox{\rm Exc} (f, \p_{1}))^{-1} $.
These integers are
\begin{equation} \label{k2}
k = e_1 - s  + [ {c_1}r/ {n_1}  ], \quad  {i_0} = k m_1 - r d_1 ,
\quad \mbox{ and } \quad {i_1} = c_1 r -  n_1  [ {c_1} r/ {n_1} ].
\end{equation}
In particular,  ${i_1}=0$ if and only if $r= p n_1$ for some
integer $p$.
\end{lemma}
{\em Proof}. We use that $\mbox{\rm Exc} (f, \p_{1})  = y_1^{n}
\circ \p_1$ by (\ref{EJ}) and that $f_2^s \circ \p_1 = (y_1^{s
n_1} \circ \p_1 ) {y_2}^s$. Hence we deduce that $  \mbox{\rm Exc}
(f, \p_{1})  {y_2}^s = ( y_1 ^{n - s n_1} f_2^s) \circ \p_1$.
Since $s< e_1$ we have that $n -s n_1= n_1 (e_1 - s)$. Then we
apply Lemma \ref{dos} for $r \geq 0$ and $l = e_1 - s >0$.
 $ \Box$

{\em Proof of Proposition \ref{comb2}.} We prove the result by
induction on the number $g$ of monomial maps in the local toric
embedded resolution, with respect to the line $\{ x_1 =0\}$. The
case $g=1$ is proved in Lemma \ref{k1}.
By induction using
(\ref{EJ}), we have that if $(r, s) \in \Z^2_{\geq 0} $ and if $s<
e_{j-1}$ there exist unique integers $k_{3}, \dots, k_j, i_0',
i_2, \dots i_j$ with $ 0 < i_0'$, $ 0 \leq i_2 < n_2$, $ \dots$, $
0 \leq i_{j-1} < n_j$ , $i_j = s$ such that
\begin{equation} \label{M}
 u_{3}^{k_{3}} \cdots u_j^{k_j} x_{j}^{r} y_{j}^{s}
=  (( x_2 ^{i_0'} y_2^{i_2} (f_3^{(1)})^{i_3} \cdots
  (f_j^{(1)})^{i_{j}} ) \circ \p_2 \cdots \circ \p_{j-1}) \,
 ({\mbox{\rm Exc} (f^{(2)}, \p_{2} \circ \dots \circ \p_{j-1} )} )^{-1}.
\end{equation}

We show that there exist  unique integers $0 < k_2, i_0$ and $0
\leq i_1 < n_1$ such that
\begin{equation} \label{m-rel}
 u_2^{k_2}  x_2 ^{i_0'} y_2^{i_2}
(f_3^{(2)})^{i_3} \cdots
  (f_j^{(2)})^{i_{j}}  = ((x_1^{i_0} \,
y_1^{i_1}  \, f_2^{i_2} \cdots f_j^{i_j}) \circ \pi_1 ) \,  (\mbox{\rm Exc} (f^{(1)}  , \p_{1} ))^{-1},
\end{equation}
By (\ref{EJ}) we have that: $ y_2^{i_2} (f_3^{(2)})^{i_3} \cdots
(f_j^{(2)})^{i_{j}} = ((y_1^{q} f_2^{i_2} \cdots f_j^{i_j}) \circ
\p_1 ) \,   (\mbox{\rm Exc} (f^{(1)}, \p_{1} ))^{-1} $,  where the
integer
\begin{equation} \label{m-qu}
q := n_1 ( e_1 - i_2 - n_2 i_3 - n_2 \cdots n_{j-1} i_j) = n_1
\left( n_2  ( \cdots ( n_{j-1}  (  e_{j-1} - i_{j} ) -i_{j-1} )
\cdots )  - i_1 \right)
\end{equation}
is a positive multiple of $n_1$ by the inequalities (\ref{lerdos}).
Then we apply Lemma \ref{dos}. $\Box$

\begin{remark}      \label{llidos}
Given the integer $e_{j-1}$ and the pairs  $(n_1, m_1), \dots, (n_{j-1}, m_{j-1})$ then
a pair $(r,s)$ with $r \geq 0$ and $s < e_{j-1}$, and the integers (\ref{lerdos}) such that
(\ref{ler}) holds, determine each other by Lemma \ref{lift} and the proof of Proposition \ref{comb2}.
\end{remark}

\begin{definition} \label{opq}  If $0 \leq r $ and if $0 \leq s < e_{j-1}$ we define
a monomial in $x, f_1, \dots, f_j$ by:
\begin{equation} \label{mjrs}
{\mathcal{M}}_j (r,s) := x^{i_0} f_1^{i_1} \cdots f_{j}^{i_j}
\end{equation}
by relation (\ref{ler}) in Proposition \ref{comb2}. We  use the
notation $\mathcal{M}_1 (r,s)$ for $x_1^r y_1^s$.  We denote the
term  $f_j^{e_{j-1}}$ by $\mathcal{M}_j (0, e_{j-1})$. We denote
the term ${\mathcal{M}}_j (r,s)$ by ${\mathcal{M}}_{j, f} (r,s)$
to emphasize the dependency with the series $f(x_1, y_1)$ defining
the plane branch $(C, 0)$.
\end{definition}
\begin{example} \label{11} The following table indicates some terms $\mathcal{M}_2(r,s)$  in the case of Example \ref{10}.
\begin{center}
\begin{tabular}{|c|c|c|c|c|c|}
\hline
  $(r,s)$ &   $(0,0) $ & $(0,1)$ & $(0,2) $& $(1,1)$ &$ (1,0)$
\\
\hline $ \mathcal{M}_2 (r,s)$  & $x_1^9$ & $x_1^6 f_2 $& $x_1^3
f_2^2 $& $x_1^5 y_1 f_2$ & $x_1^8 y_1$
\\
\hline
\end{tabular}
\end{center}
For instance, we have that $\mathcal{M}_2 (1, 1)=  x_1^5 y_1 f_2$,
since $x_1^5 y_1 f_2 \circ \p_1 = \mbox{\rm Exc} (F^{(1)}, \p_1)
u_2^2 x_2 y_2$, where $\mbox{\rm Exc} (F^{(1)} , \p_1) =  u_2 ^6
x_2 ^{18 }$ by Example \ref{10}. Notice also that the analytic
function $x_2 y_2 \mbox{\rm Exc} (F, \p_1)$ on $Z_2$ is equal to $
( x_1^{-1} y_1^5 f_2) \circ \p_1$, i.e., it is the transform by
$\p_1$ of a {\em meromorphic} function. Both of the following
formulas
\[ y_1^6 \circ \p_1 = \mbox{\rm Exc} (F^{(1)}, \p_1) \mbox{ and  }
x_1^9 \circ \p_1 = \mbox{\rm Exc} (F^{(1)}, \p_1) u_2^3\]  seem to
correspond to (\ref{ler}) in the case $(r,s) =(0,0)$, however the
term $y_1^6$ is not of the form prescribed by the inequalities
(\ref{lerdos}), hence the first formula is not the one considered
by Lemma \ref{comb2}.

\begin{lemma} \label{m-min} If $0 \leq r$ and $s < e_{j-1}$, we have
that:
\[ ({\mathcal{M}}_j (r,s), f )_0  = e_{j-2} \bar{b}_{j-1} + r e_{j-1} + s \,(
{\bar{b}}_{j+1} - n_j {\bar{b}}_j), \quad \mbox{ for } j=2,\dots, g+1. \]
\end{lemma}
{\em Proof.} By Lemma \ref{comb2}   we have that:
\[{\mathcal{M}}_j (r,s)  \circ \p_1 \circ \cdots \circ \p_{j-1} =
\mbox{\rm Exc} (f, \p_{1} \circ  \cdots \circ \p_{j-1} ) \,
u_2^{k_2} \cdots u_j^{k_j} x_j^{r} y_j^{s}.\] By Proposition
\ref{mod2} and Lemma \ref{j-semi} we deduce that: \[ (
{\mathcal{M}}_j (r,s) , f)_0 = \ord_\t (\mbox{\rm Exc} (f, \p_{1}
\circ  \cdots \circ \p_{j-1} )) +     r e_{j-1} + s \,(
{\bar{b}}_{j+1} - n_j {\bar{b}}_j) \quad \Box.\]
 \end{example}

\begin{lemma} \label{newton} If $0 \leq r$ and $0 \leq s < e_{j-1}$
the Newton polygon of a term ${\mathcal{M}}_j (r,s)$, with respect
to the coordinates $(x_1, y_1)$, is contained in the Newton
polygon ${\mathcal N} (f(x_1, y_1))$ for $2 \leq j \leq g+1 $. It
is contained in the interior of ${\mathcal N}(f(x_1, y_1))$ unless
$j=2$, $r=0$ and $ 0 \leq s < e_1$.
\end{lemma}
{\em Proof.}  If $j=2$ we have that ${\mathcal{M}}_2 (r,s ) =
x_1^{i_0} y_1^{i_1} f_2^{s}$ by Lemma \ref{k1}. The vector
$\vec{v}:= ( i_0+ s m_1 , i_1)$ is a vertex of the Newton polygon
of ${\mathcal{M}}_2(r,s)$ and $\vec{w} := (\bar{b}_1,0)$ is a
vertex of the only compact edge $\Gamma_1$ of  ${\mathcal N}
(f(x_1, y_1))$. Notice that if $s =0$ then Newton polygon of
${\mathcal{M}}_2(r,s)$ has only one compact face $\{ \vec{v} \}$,
otherwise it has only one compact edge which is parallel to
$\Gamma_1$ (see Subsection \ref{toric-res}). The vector $\vec{p}_1
= (n_1, m_1)$ is orthogonal to $\Gamma_1$ hence we deduce the
inequality:
  \[
n_1 \bar{b}_1 = \langle \vec{p}_1, \vec{w} \rangle \leq \langle
\vec{p}_1, \vec{v} \rangle = n_1 i_0 + s  n_1 m_1 + m_1 i_1 = e_1
n_1 m_1 + r(m_1 c_1 - n_1 d_1) \stackrel{(\ref{ave-i})}{=}
n_1\bar{b}_1 + r,
 \]
 using (\ref{k2}). Equality holds in formula above if and only if
$r=0$.

If $j >2$ we follow the proof of Proposition \ref{comb2}: there
exist integers $0< i_1 \leq n_1$, $0 < i_0',  k_2 $ such that
(\ref{m-rel}) holds. By Lemma \ref{dos} we have that $k_2 = l +
[\frac{c_1 i_0'}{n_1}]$ where the integer $l$ is $l := e_1 - i_2 -
n_2 i_3 - \cdots - n_2 \cdots n_{j-1} i_j$. The vector $\vec{v}:=
(i_0 + m_1 (e_1 -l) ,i_1)$ is a vertex of ${\mathcal N}
({\mathcal{M}}_j(r,s))$. By  the construction the Newton polygon
of ${\mathcal{M}}_j(r,s)$ has at most one compact edge, which is
in addition  parallel to $\Gamma_1$. We deduce from a simple
calculation using (\ref{m-k3}) that:
\begin{equation} \label{m-be}
n_1 \bar{b}_1 \leq \langle \vec{p}_1, \vec{v} \rangle = n_1
\bar{b}_1 + i_0'.
\end{equation}
By the proof of Proposition \ref{comb2} we  have that $i_0' >0$,
hence the inequality (\ref{m-be}) is strict. $\Box$

\begin{remark}  \label{newton-rem} By induction using the same arguments as in Lemma \ref{newton} we check that
if $1 \leq i < j$,  $0  \leq r$, and $0 \leq s  < e_{j-1}$ that the Newton polygon of
$({\mathcal{M}}_ j(r,s) \circ \p_1\circ \cdots \circ \p_{i-1})
\, ( \mbox{\rm Exc } (f, \p_1\circ \cdots \circ \p_{i-1} )) ^{-1} $
with respect to the coordinates $(x_i, y_i)$, is contained in $\mathcal{N} (f^{(i)})$.
  It is contained in the interior of
${\mathcal N}(f^{(i)})$ unless $j=i+1$, $r=0$ and $ 0 \leq s < e_{i}$.
\end{remark}

\section{Irreducibility and equisingularity criterions} \label{k-exp}

Abhyankar's irreducibility criterion gives an affirmative answer
to a question of Kuo mentioned in \cite{Abhyankar}: "Can we decide
the irreducibility of a power series $F(x, y)$ without blowing up
and without getting into fractional power series ?" We have
revisited the Abhyankar's criterion in the light of toric geometry
methods. In particular, our proof explains that if $F$ is
irreducible,  some information on the Newton polygons of the strict
transform of $F$ at the infinitely near points of the toric
resolution can be read from the expansions in certain semi-roots
of $F$. See \cite{Cossart} and \cite{Cossart-Saskatoon}, for an
extension of this criterion to the case of base field of positive
characteristic.       As an application we give an equisingularity
criterion for an equimultiple family of plane curves to be
equisingular to a plane
 branch (See Section \ref{equi-crit}).

\subsection{Straight line conditions in the toric resolution}  \label{st-line}

We consider a plane branch $(C,0)$ together with its local toric
resolution. We keep notations of Section \ref{toric-res} (see also
Notation \ref{not-semi}).  We give some precisions on the $( f_1,
\cdots, f_j)$-expansion of $f$ (see Proposition \ref{g-expand}).
We have that the $(f_1, \cdots, f_j)$-expansion of $f$ is of the
form:
\begin{equation} \label{expansion-j}
f = f_j^{e_{j-1}} + \sum_{I=(i_1, \dots, i_{j})} \a_{I} (x_1)
f_1^{i_1} \cdots f_{j}^{i_{j}}, \mbox{ with } \a_I (x_1) \in \C \{
x_1 \},
\end{equation}
with $ 0 \leq i_1 < n_1$, $\dots$, $ 0\leq i_{j-1} < n_{j-1}$, $0
\leq i_j < e_{j-1}$,  for $2\leq j \leq g$.

By expanding the coefficients of the terms in (\ref{expansion-j}),
as series in $x_1$,  we obtain the following expansion
\begin{equation} \label{x-expansion-j}
f= f_j^{e_{j-1}} + \sum_{J = (i_0, \dots, i_j) } \b_J \, x_1^{i_0}
\, f_1^{i_1} \cdots f_j^{i_j} \mbox{ with } \b_J \in \C,
\end{equation}
which we call the {\em $(x_1, f_1, \dots, f_j)$-expansion} of $f$.
The main result of this section is the following (see Definition
\ref{opq}).

\begin{theorem} \label{straight}
The $(x, f_1, \dots , f_j)$-expansion of $f$, for $j=2, \dots, g$,
is of the form: \[ f=  f_j^{e_{j-1}}  + \sum_{(r, s)}  c_{r,s}  \,
{\mathcal{M}}_j (r,s),\]  where $c_{r,s} \in \C$ and  the pairs
$(r, s)\in \Z^2$ verify that \[ 0 < r, \quad 0 \leq s< e_{j-1},
\quad e_{j-1} (\bar{b}_{j } - n_{j-1} \bar{b}_{j-1}) \leq  r
e_{j-1} + s (\bar{b}_{j} - n_{j-1} \bar{b}_{j-1})  .\] Among the
terms of this expansion with minimal intersection multiplicity
with $f$ there exist $f_j^{e_{j-1}}$ and ${\mathcal{M}}_j (
\bar{b}_{j +1} - n_j \bar{b}_j , 0)$. Moreover, if $j= g-1$ these
two terms are exactly the terms with minimal intersection
multiplicity with $f$.
\end{theorem}

Before entering into the proof of Theorem \ref{straight} we
discuss the following Propositions.
\begin{proposition} \label{newton-j}
 If $j > 1$ and the coefficient $\a_I (x_1)$ in (\ref{expansion-j}) does
 not vanish then the Newton polygon of the term $\a_{I} (x_1) f_1^{i_1} \cdots
f_{j}^{i_{j}}$ is contained in the interior of the Newton polygon
of $f$.
\end{proposition}
{\em Proof.} Since $\deg f = e_{j-1} \deg f_j$ and both are monic
polynomials we have that the term $f_j^{e_{j-1}}$ appears in the
$(f_1, \cdots, f_j)$-expansion of $f$ with coefficient one.

For an index $I= (i_1, \dots, i_j)$ appearing in
(\ref{expansion-j}) we denote by ${\mathcal{M}}_I $ the term $\a_I
(x) f_1^{i_1} \cdots f_{j}^{i_{j}}$. By  Remark \ref{one-edge},
the  Newton polygon  $\mathcal{N} ({\mathcal{M}}_I) $ of
${\mathcal{M}}_I$ has only one compact face $\Gamma_I$ of maximal
dimension which is parallel to the compact face $\Gamma_1$ of
${\mathcal{N}}(f(x_1, y_1))$. The vector $\vec{p_1} = (n_1, m_1)$,
which was defined in Section \ref{toric-res}, is orthogonal
$\Gamma_1$.

We set also the numbers  \[ B_I := \min \{ \langle \vec{p_1},
\vec{u} \rangle \mid \vec{u} \in \mathcal{N} ({\mathcal{M}}_I) \}
\mbox{ and }  q_I := \ord_{y_1} { ( f_1^{i_1} \cdots f_{j}^{i_{j}}
) }_{|x_1 =0}, \]  for $I$ appearing in  the expansion
(\ref{expansion-j}) with non-zero coefficient. The numbers $q_I$
defined above, are all distinct by Remark \ref{injective} applied
to  $0 \leq i_1 < n_1$, $\dots$, $ 0\leq i_{j-1} < n_{j-1}$ and $0
\leq i_j < e_{j-1}$.

Suppose that there exists an index $\tilde{I} = (\tilde{\imath}_1,
\dots, \tilde{\imath}_j)$   with $\a_{\tilde{I}} \ne 0$, such that
the polygon $\mathcal{N} ({\mathcal{M}}_{\tilde{I}}) $ is not
contained in ${\mathcal{N}}(f(x_1, y_1))$. This holds if and only
if $B_{\tilde{I}} < \, \min \{ \langle \vec{p_1}, \vec{u} \rangle
\mid \vec{u} \in \mathcal{N} (f) \}$. Hence
$\mathcal{M}_{\tilde{I}} $ is not equal to $f_j^{e_{j-1}}$, since
$\mathcal{N}(f) = \mathcal{N} (f_j^{e_{j-1}})$ by Remark
\ref{one-edge}.
We can suppose in addition that $B_{\tilde{I}}$ is
the minimal number of this form. Moreover, we can assume that
$\tilde{I}$ has the following property: if the index $I = (i_1',
\dots, i_j') \ne \tilde{I}$, which appears in (\ref{expansion-j})
with non zero coefficient, verifies that $B_{\tilde{I}} = B_{I}$
then $q_{\tilde{I}} > q_{I}$.
 If
$(r,s) \in \Gamma_{\tilde{I}} \cap \Z^2$, the sum $K_{r,s}$ of the
coefficients of the term $x^r y^s$ in $\a_{{I}}
{\mathcal{M}}_{I}$, for those indices $I$ with $B_I =
B_{\tilde{I}}$, must vanish. But if $(r,s)$ is the vertex of
$\Gamma_{\tilde{I}}$ with $s = q_{\tilde{I}}$ then we obtain that
$K_{r,s} $ is the initial coefficient of the series
$\a_{\tilde{I}}$, a contradiction. Thus, for all index $I$
appearing in  (\ref{expansion-j}) we have the inclusion
$\mathcal{N} ({\mathcal{M}}_{{I}}) \subset  \mathcal{N} (f (x_1,
y_1) ) $.

By Remark \ref{one-edge}  the symbolic restrictions of $f$ and of
$f_j^{e_{j-1}}$, to the compact face $\Gamma_1$ of the Newton
polygon coincide. Suppose that  there exists an index $I$
appearing in the expansion  (\ref{expansion-j}) with non zero
coefficient such that ${\mathcal{M}}_{I} \ne f_j^{e_{j-1}}$ and
$B_I = \min \{ \langle \vec{p_1}, u \rangle \mid u \in \mathcal{N}
(f) \}$. In this case for any $(r,s) \in \Gamma_1 \cap \Z^2$ the
sum of the coefficients of the terms $x^r y^s$ in
${\mathcal{M}}_{I'}$, for those $I'$ with $B_I = B_{I'}$ and
${\mathcal{M}}_{I'} \ne  f_j^{e_{j-1}}$, must vanish. We argue as
in the previous case to prove that this cannot happen. $\Box$

\begin{lemma} \label{num}
Let $\mathcal{M}_J = x_1^{i_0} f_1^{i_1} \cdots f_{j}^{i_{j}} $ be
a term in the  in the expansion (\ref{x-expansion-j}) with
non-zero coefficient corresponding to the index $J = (i_0, \dots,
i_j)$. Set $q_J := i_1 + n_1 i_2 + \cdots + n_1 \cdots n_{j-1}
i_j$. We can factor ${\mathcal{M}}_J \circ \p_1$ as:
\begin{equation} \label{uss}
({\mathcal{M}}_J \circ \p_1 ) \, (\mbox{\rm Exc} (f, \p_1) )^{-1}=
\, u_2 ^{k_2 (J) } \, x_2^{i_0' (J)} \,  (f_2^{(2)})^{i_2} \cdots
(f_j^{(2)})^{i_j},
\end{equation}
where $i_0' (J) = n_1 i_0 - m_1 (e_0 - q_J)  >  0$ and $k_2 (J)  =
c_1 i_0 - d_1 (e_0 - q_J) >  0$.
\end{lemma}
{\em Proof.} Notice that $q_J$ is the degree in $y$ of the term
${\mathcal{M}}_J $. By Proposition \ref{newton-j} the Newton
polygon of the term ${\mathcal{M}}_J$ is contained in the interior
of the Newton polygon of $f$. This implies that  $\vec{v}_J =
(i_0, q_J)$ is a vertex of the Newton polygon of $
{\mathcal{M}}_J$  and $\langle \vec{p}_1, \vec{v}_J \rangle
> e_0 m_1$. This implies that  $i_0' (J) >0$. We deduce from
this that $k_2(J) >0$ and (\ref{uss}) holds. $\Box$

 We
obtain the following expansion from (\ref{x-expansion-j}), by
factoring out $\mbox{\rm Exc} (f, \p_1)$ from $f \circ \p_1$:
\begin{equation} \label{2x-expansion-j}
f^{(2)} = (f_j^{(2)})^{e_{j-1}} + \sum_{J = (i_0, \dots, i_j) }
c_J \, u_2 ^{k_2 (J) } \, x_2^{i_0' (J)} \,  (f_2^{(2)})^{i_2}
\cdots (f_j^{(2)})^{i_j}.
\end{equation}

The following expansion is obtained from (\ref{2x-expansion-j}) by
collecting the terms with the same index $I' = (i_2, \dots, i_j)$:
\begin{equation} \label{j-2}
f^{(2)} = (f_j^{(2)})^{e_{j-1}} + \sum_{I'=(i_2, \dots, i_j)}
\a_{I'}^{(2)} (x_2, u_2) \, ( f_2^{(2)} ) ^{i_2} \cdots
(f_{j}^{(2)}) ^{i_{j}}.
\end{equation}
By (\ref{xi-rel}) the coefficient $\a_{I'}^{(2)} (x_2, u_2)$,
viewed in $ \C \{x_2, y_2 \}$, is of the form:
\begin{equation} \label{coefficient-2}
\a_{I'}^{(2)}= \epsilon^{(2)}_{I'} x_2^{r_2(I')} \mbox{ with }
r_2(I') > 0  \mbox{ and } \epsilon_{I'}^{(2)} \mbox{ a unit in }
\C \{x_2, y_2 \}.
\end{equation}

\begin{definition}
We call the expansion (\ref{2x-expansion-j}) (respectively
(\ref{j-2})) the $(u_2, x_2, f_2^{(2)}, \dots,
f_j^{(2)})$-expansion (respectively $(f_2^{(2)}, \dots,
f_j^{(2)})$-expansion) of $f^{(2)}$.
\end{definition}

\begin{proposition}  \label{newton-j-2} Suppose that $2 \leq j \leq g$.
Let us consider an index $I' = (i_2, \dots, i_j)$ appearing in the
expansion (\ref{j-2}) with coefficient $\a_{I'}^{(2)} (x_2, u_2)
\ne 0$.  Denote by $q_{I'}$ the order in $y_2$ of the series
 $( ( f_2^{(2)} ) ^{i_2} \cdots (f_{j}^{(2)}) ^{i_{j}}) _{|x_2
 =0}$.
\begin{enumerate}
\item[(i)] For any pair $I'_1\ne I'_2$ of indexes  in (\ref{j-2})
with $\a_{I_1'}^{(2)} \a_{I'_2}^{(2)} \ne 0$ we have that
$q_{I'_1} \ne q_{I_2'}$.

\item[(ii)] If $j >2$ and $\a_{I'}^{(2)} (x_2, u_2) \ne 0$ the
Newton polygon of the term $ \a_{I'}^{(2)} (x_2, u_2)  ( f_2^{(2)}
) ^{i_2} \cdots (f_{j}^{(2)}) ^{i_{j}})   $ (with respect to the
coordinates $(x_2, y_2)$) is contained in the interior of
$\mathcal{N} (f^{(2)} (x_2, y_2))$.

\end{enumerate}
\end{proposition}
{\em Proof.} The assertion on the orders  in $y_2$ of the series
$(( f_2^{(2)} ) ^{i_2} \cdots (f_{j}^{(2)}) ^{i_{j}} )_{|x_2 =0}$
is consequence of Remark \ref{injective} with respect to the
integers $n_2, \dots, n_{j-1} > 1$.

For the second assertion notice that the Newton polygon with
respect to the coordinates $(x_2, y_2)$ of a term $
{\mathcal{M}}_{I'}^{(2)} := \a_{I'}^{(2)} (x_2, u_2) \, (
f_2^{(2)} ) ^{i_2} \cdots  ( f_j^{(2)} ) ^{i_j}$, appearing in the
expansion (\ref{j-2}), has at most one compact face which is
parallel to $\Gamma_2$. We deduce that the Newton polygon of
${\mathcal{M}}_{I'}^{(2)}$ is contained in the interior of
$\mathcal{N} (f^{(2)} (x_2, y_2))$ by repeating the argument of
Proposition \ref{newton-j} combined with Remark \ref{one-edge}
(ii). $\Box$

{\em Proof of Theorem \ref{straight}.} Let ${\mathcal{M}}_I :=
x^{i_0} f_1^{i_1} \cdots f_j^{i_j}$ be a monomial appearing in
(\ref{x-expansion-j}). Using Proposition \ref{newton-j},
Proposition \ref{newton-j-2} and Lemma  \ref{num} we deduce
inductively that $ ({\mathcal{M}}_I \circ \p_1 \circ \cdots \circ
\p_{j-1}) \, (\mbox{\rm Exc} (f, \p_1 \circ \cdots \circ
\p_{j-1}))^{-1}= u_2^{k_2(I) } \cdots u_j ^{k_j (I)} x_j^{r(I)}
y_j^{s(I)}$, where $k_2(I), \dots, k_j(I), r(I) > 0$ and $s(I) =
i_j$. It follows that we have an expansion:
\begin{equation} \label{x-expansion-jj}
f^{(j)}  = y_j^{e_{j-1}} + \sum_I c_{I} \, u_2^{k_2(I)} \cdots u_j
^{k_j(I)} \,  x_j^{r(I)} \, y_j^{s(I)}.
\end{equation}
 By the unicity statement in Proposition \ref{comb2} it
follows that $\mathcal{M}_I = \mathcal{M} (r(I), s(I))$, hence, if
$I \ne I'$ are two different indices appearing in
(\ref{x-expansion-j}), then $(r(I), s(I)) \ne (r(I'), s(I') )$. By
(\ref{xi-rel}) the term $u_2^{k_2(I)} \cdots u_j ^{k_j (I)}$ is a
unit viewed in $\C \{ x_j, y_j \}$,  therefore the Newton polygon
of $f^{(j)} (x_j, y_j)$ is equal to the convex hull of the set,
$\cup_{I} (r(I) ,s(I)) + \R_{\geq 0}^2 $. By Proposition
\ref{mod2} this polygon has vertices $(0, e_{j-1})$ and
$(\bar{b}_{j} - n_{j-1} \bar{b}_{j-1} , 0)$. If $j =g$, these two
vertices are the unique integral points in the Newton polygon. By
Lemma \ref{m-min}, the exponents $(r,s) \in \Gamma_j$ correspond
to terms ${\mathcal{M}}_j (r,s)$ with minimal intersection
multiplicity with $f$ at the origin. $\Box$

We deduce from Theorem \ref{straight}  the following Corollary,
where the coefficient $\theta_j$ is the same as the one appearing
in Formula (\ref{fg}) at level $j$.
\begin{corollary} \label{coro}
If  $j \in \{ 2, \dots, g \} $, the $(x, f_1, \dots ,
f_j)$-expansion of $f_{j+1}$ is of the form (cf. Definition
\ref{opq})
\begin{equation} \label{mas}  f_{j+1} = f_j^{n_j}  -
\theta_{j}
 \mathcal{M}_{j, f_{j+1}} (m_j, 0) +
 \sum_{(r, s)} c_{r,s} \, {\mathcal{M}}_{j, f_{j+1}}  (r,s),
 \end{equation} where
$(r, s)$ above verify that $ 0 < r$, $0 \leq s < n_j$ and $ n_j
m_j < n_j r + m_j s$.
\end{corollary}

\begin{remark} \label{normalized}
In some cases it may be useful to have $\theta_1 = \dots =
\theta_g = 1$. We can reduce to this case by replacing the terms
$(x_1, f_1, \dots, f_{g})$ by $(\eta_0 x_1, \eta_1 f_1, \dots,
\eta_g f_g)$ for some suitable constants $\eta_0, \dots, \eta_g
\in \C^*$.
\end{remark}
To see this, by a change of coordinates of this form, we can assume that the  image
of $x, f_1, \dots, f_g$  in the integral closure $\C \{ \t \} $ of the
algebra of $(C,0)$, are series with
constant term equal to one. By Lemma \ref{m-min} we have that
$\bar{b}_{j+1} = \ord_\t f_{j+1} (x_1 (\t), y_1 (\t)) > n_j
\bar{b}_j = \ord_\t ( f_{j}^{n_j} (x_1 (\t), y_1 (\t)) = \ord_\t (
\mathcal{M}_{j, f_{j+1}} (m_j, 0) ) $ and  $ n_j \bar{b}_j <
\ord_\t ( \mathcal{M}_{j, f_{j+1}} (r, s) )$, for those pairs $(r,
s)$ appearing in  (\ref{mas}). We deduce  by a standard valuative
argument that in this case $\theta_j= 1$.

\subsection{Abyankar's generalized Newton polygons,
straight line condition and irreducibility criterion}

We follow the presentation given by Assi and Barile in \cite{Assi}
of results in \cite{Abhyankar}.

\subsubsection{Generalized Newton polygons}

Given a sequence  $\underline{\bar{B}}:= (\bar{B_0} , \bar{B}_1,
\dots, \bar{B}_G)$  of positive integers with $\bar{B}_1 < \cdots <
\bar{B}_G$, we associate to them sequences $E_{j} = \mbox{\rm gcd
} ( \bar{B_0}, \bar{B}_1, \dots, \bar{B}_j)$ and $N_0 =1$, $N_j =
E_{j-1} / E_j$, for $j=0, \dots, G$. Notice that if $
\underline{\bar{B}}$  a characteristic sequence of generators of
the semigroup $\Lambda_C$ associated to a plane branch $(C, 0)$,
we set $g = G$ and we have with the Notations of the first section
that  $E_j = e_j$ and $N_j = n_j$, for $j=0, \dots, g$.

Let $F$ be a Weierstrass polynomial of the form:
\begin{equation}    \label{candidato}
         F = y^{N} + \sum_{i=2}^N A_i (x) y^{N-i} \in \C \{ x \} [y]
\end{equation}
We assume that $y$ is an approximate root of $F$ since the
coefficient of $y^{N-1}$ is equal to zero. We denote by     $F_j$
the approximate root of $F$ of degree $N_0 \dots N_{j-1}$, and by
$\underline{F}_j$ the sequence $(F_1, \dots, F_{j})$ for $j
=1,\dots, G+1$ and $\underline{F} = \underline{F}_{G+1}$.

Let $P \in \C \{ x \} [y] $ be a monic polynomial. The $(F_1,
\dots, F_{G+1})$-expansion of $P$ is of the form  $P = \sum_I \a_I
(x) {F}_1^{i_1} \cdots {F}_{G}^{i_{G}} {F}_{G+1}^{i_{G+1}} $ (see
Proposition \ref{g-expand}).
 The {\em formal intersection
multiplicity} of $P$ and $\underline{F}$, with respect to the
sequence $\underline{\bar{B}}$ is defined as
\begin{equation}
\mbox{\rm formal }_{\underline{\bar{B}}} (  P,  \underline{F} )
 := \min \{ \sum_{j=0} ^G i_j \bar{B}_j      \mid I= (i_1,\dots, i_{G}, 0) ,
 \a_I (x) \ne 0 \}.
\end{equation}
Notice that when this value is $< + \infty$, it is reached at only
one coefficient.

Let $P, Q \in \C \{ x \} [y] $ be two monic polynomials of degrees
$p, q$ with $p =mq$. We have the $Q$-adic expansion of $P$ is of
the form: $
          P = Q^m + \a_1 Q^{m-1} + \cdots +  \a_m$.
The {\em generalized Newton polygon} $ \mathcal{N} (P, Q,
\underline{\bar{B}},     \underline{F} )$ of $P$ with respect to $Q$ and the sequences
$\underline{\bar{B}}$
and $\underline{\bar{F}}$ is the convex hull of the set:
\begin{equation} \label{gen-newton-pol}
\bigcup_{k=0}^m       ( \mbox{\rm formal }_{\underline{\bar{B}}} (  \a_k ,
\underline{F} ), (d-k) \mbox{\rm formal }_{\underline{\bar{B}}} (
Q ,  \underline{F} ))
   + \R^2_{\geq 0}.
\end{equation}

\subsubsection{Abhyankar's irreducibility criterion}

To a monic polynomial $F$ of the form (\ref{candidato}) it is
associated a sequence $\bar{B}$  as follows: Set $ \bar{B}_0 = E_0
:= N$, $F_1 = y$, $\bar{B}_1 = (F_1, F)_0$, $E_1 = \mbox{ gcd }
(\bar{B}_0, \bar{B}_1) $  and $ N_1 := E_0 / E_1$. Then, for $j
\geq 2$  the integers  $E_0 \dots, E_j =     \mbox{ gcd }
(\bar{B}_0, \dots, \bar{B}_{j})$ and $N_1, \dots, N_{j-1}$ are
defined by induction. We set $B_{j+1} = (F, F_{j})_0$, where
$F_{j}$ denotes  the approximate root of $F$ of degree $N_1 \cdots
N_{j-1}$.

\begin{theorem}{\rm (\cite{Abhyankar})} \label{assi}
With the above notations  the polynomial $F \in \C \{ x\} [y]$ is
irreducible if and only if the following conditions hold:
\begin{enumerate}
\item[(i)] there exists an integer $G \in \Z_{> 0}$ such that
$E_{G} =1$,

\item[ii)] $\bar{B}_{j+1} > N_j \bar{B}_j$ for $j=1, \dots G-1$,

\item[(iii)]  {\em (straight line condition)} the generalized
Newton polygon
          $\mathcal{N} (F_{j+1}, F_j,   \frac{1}{E_j} \underline{\bar{B}}_j,     \underline{F}_j ) $
has only one compact edge with vertices
$(\frac{1}{E}_j N_j \bar{B}_j, 0)$ and   $(0, \frac{1}{E}_j N_j \bar{B}_j)$.
\end{enumerate}
\end{theorem}
{\em Proof.} We prove first that if $F$ verifies the conditions of
the theorem  then $F$ is irreducible. By the straight line
condition the vertices of the generalized Newton polygon,
$\mathcal{N} (F_{j+1}, F_j,   \frac{1}{E_j} \underline{\bar{B}}_j,
\underline{F}_j )$,  correspond to the terms $F_j^{n_j}$ and
$\a_0^{(j)} (x)  F_1^{ \eta_{1}^{(j)} } \cdots F_{j-1}^{
\eta_{j-1}^{(j)} }$  of the $(F_1,\dots, F_j)$ expansion of
$F_{j+1}$, where $\ord_x \a_0^{(j)} (x) = \eta_{0}^{(0)}  \geq 0$
and $0 \leq \eta_i^{(j)}  < N_i$ for $i =1, \dots, j-1$. The
straight line condition implies that $\frac{1}{E_j} N_j \bar{B}_j
= \frac{1}{E_j} (\eta_0^{(j)} \bar{B}_0 + \cdots +
\eta_{j-1}^{(j)} \bar{B}_{j-1}) $. It follows that $N_j \bar{B}_j
$ belongs to the semigroup generated by $\bar{B}_0, \dots,
\bar{B}_{j-1}$. By Lemma \ref{tee} this numerical condition
together with  (i) and (ii) guarantee that the semigroup generated
by  $\bar{B}_0, \dots, \bar{B}_{G}$ is the semigroup of a plane
branch.

If $\alpha(x) F_1^{\eta_1} \cdots F_{j-1}^{\eta_{j-1}}$  with $\mathrm{ord}_x \alpha(x) =\eta_0 \geq 0$, 
is a term appearing in the 
$(F_1, \dots, F_j)$-expansions of $F_{j+1}$, different from $F_j^{N_j}$ and 
$\alpha_0^{(j)} (x) F_1^{\eta_1 ^{(j)}} \cdots F_{j-1}^{\eta_{j-1} ^{(j)}}$, then 
the straight line condition implies that 
\begin{equation}  \label{**}
 (1/E_j)  N_j \bar{B}_j < (1/E_j) ( \eta_0 \bar{B}_0 + \cdots + \eta_{j-1} \bar{B}_{j-1} ).
\end{equation}
The inequality  (\ref{**}) is strict since $\sum_{i=0}^{j-1}   \eta_i \bar{B}_i \ne 
\sum_{i=0}^{j-1}   \eta_i^{(j)} \bar{B}_i$ by the numerical properties of the 
semigroup $\Gamma$ generated by $\bar{B}_0, \dots, \bar{B}_G$, see Lemma 1.15. 

\medskip

We re-embed the germ  $(C, 0)$ defined by $F =0$, in $(\C^{G+1},0)$ by setting 
\begin{equation} \label{*}
u_0 = x, \quad u_1 = F_1, \quad \dots \quad, u_G = F_G. 
\end{equation}
We also set the weight of $u_i$ equal to $\bar{B}_i$, for $i =0, \dots, G$. 
The equations defining the embedding of $(C, 0)$ are obtained by making the replacement (\ref{*}) 
in the $(F_1, \dots, F_j)$-expansion of $F_{j+1}$ for $j=1, \dots, G$.  
The inequalities of the form (\ref{**}) together with 
$N_j \bar{B}_j =  \sum_{i=0}^{j-1}   \eta_i^{(j)} \bar{B}_i  $ for $j =1, \dots, G$ 
and  $ N_j \bar{B}_j < \bar{B}_{j+1}$  for $j= 1, \dots, G-1$, 
are precisely the weight conditions 
on the equations of defining the embedding 
$(C, 0) \subset (\C^{G+1}, 0)$ 
indicated in Proposition 39 of \cite{GP}.
%the generic fiber of an equisingular deformation of the monomial curve $C^\Gamma$ parametrized by 
% $u_i= t^{\bar{B}_i}$, $i =0, \dots, G$ (see \cite{Rebeca, T, GP}). 
% This implies that $(C, 0)$ is an irreducible germ with semigroup $\Gamma$ 
% and hence $F$ is irreducible. 
By the proof of Theorem 6.1 in 
\cite{Rebeca} or by Theorem 2, page 1867 in \cite{GP}
one toric modication of $\C^{G+1}$  provides a simultaneous embedded resolution of both
the monomial curve $(C^\Gamma, 0)$  parametrized by 
$u_i= t^{\bar{B}_i}$, $i =0, \dots, G$ and the germ $(C, 0) \subset (\C^{G+1},0)$.  
In addition, the strict transform of both curves intersect the exceptional divisor at 
exactly one point. In particular, this implies that the normalization of $(C, 0)$ is 
smooth, and the germ $(C, 0)$ is irreducible hence $F$ is irreducible.

% Let $0 \ne F' \in  \C \{ x \} [y]$ be any polynomial of
% degree $< N = \deg F$. We consider the $(F_1,\dots, F_G)$
% expansion of $F'$:
% \begin{equation} \label{A-expand}
% F ' = \sum_I  \a_I (x) {F}_1^{i_1} \cdots {F}_{G}^{i_{G}} \mbox{ with } \a_I (x) \in \C\{x\},
% \end{equation}
% Set $i_0 = \ord_x   \a_I (x) $. The intersection multiplicities
% $(F, \a_I (x) {F}_1^{i_1} \cdots {F}_{G}^{i_{G}} )_0 =
% \sum_{j=0}^G i_j \bar{B}_j$ obtained for the different terms in
% the expansion (\ref{A-expand}) are all different (this reduces to
% an arithmetical property which can be proved similarly as Lemma
% \ref{tee}).
%  We deduce
% that $(F, F')_0 = \min  \{    \sum_{j=0}^G i_j \bar{B}_j \mid \a_I
% (x) \ne 0 \} < + \infty$. The polynomial $F$ is irreducible,
% otherwise there is an irreducible factor of $ F'$ of $F$  of
% degree $< \deg F$ and then $(F, F')_0 = +\infty $,  a
% contradiction.

Suppose now that $F$ is irreducible. Then $\bar{B}_0, \dots,
\bar{B}_G$ are the generators of the semigroup of the branch $F=0$
with respect to the line $\{ x=0 \}$. By Lemma \ref{tee}, the
first two conditions in the statement of the theorem hold
automatically. By Proposition \ref{A-M} the approximate root
$F_{j+1}$ is irreducible and define a plane branch with semigroup
generated by $\frac{1}{E_j}\bar{B}_0, \dots,
\frac{1}{E_j}\bar{B}_j$. By Theorem \ref{straight}
the Newton polygon of $F^{(j)}_{j+1} (x_j, y_j)$ has only two
vertices  $(0, N_j)$ and $(M_j, 0)$ which correspond respectively
to the terms $F_{j}^{N_j} $ and ${\mathcal{M}}_{j, F_{j+1}} (0,
M_j)$  of the $(x, F_1, \dots, F_{j})$-expansion of $F_{j+1}$ (see
Definition \ref{opq}). By Proposition \ref{mod2} and induction the
vertices of the Newton polygon of the strict transform function
  $F^{(j)} (x_j, y_j)$ are
$(\bar{B}_j - N_{j-1} \bar{B}_{j-1}, 0)$ and $(0, E_j N_j)$. It
follows that $M_j = \frac{1}{E_j}    ( \bar{B}_j - N_{j-1}
\bar{B}_{j-1} ) $. By Lemma \ref{j-semi} we have that $\ord_t
(\mbox{{\rm Exc}} (F_{j+1}, \p_1 \circ \cdots \circ \p_{j-1} )) =
\frac{E_{j-2}}{E_j} \frac{\bar{B}_{j-1}}{E_j} = \frac{1}{E_j}
N_{j-1} N_j \bar{B}_j$. By Lemmas  \ref{m-min} and \ref{j-semi} 
we deduce the equalities
\[
\begin{array}{lcl}
(F _{j+1}, {\mathcal{M}}_{j, F_{j+1}} (0, M_j))_0  & = &
\frac{1}{E_j} N_{j-1} N_j \bar{B}_{j-1} + \frac{E_{j-1}}{E_j} M_j
\\
& =  & \frac{1}{E_j}    ( N_j   N_{j-1}  \bar{B}_{j-1}      + N_j  (
\bar{B}_j - N_{j-1} \bar{B}_{j-1} )
\\
& = &
 \frac{1}{E_j}  N_j  \bar{B}_j
\\
& = &
 (F_{j+1},
     F_{j}^{N_j})_0.
\end{array}
\]
A similar computation using    Lemmas \ref{m-min} and \ref{j-semi}
proves that if
   $  {\mathcal{M}}_{j, F_{j+1}} (r, s)$  appears in the  $(x, F_1, \dots, F_{j})$-expansion of $F_{j+1}$
then    $(F_{j+1},    {\mathcal{M}}_{j, F_{j+1}} (r, s) )_0  >
\frac{1}{E_j} N_j  \bar{B}_j  $.
  $\Box$

          \begin{remark} \label{gen-aby}
We keep notations of the proof of Theorem \ref{assi}. Suppose that
$F$ is irreducible. Let $\f_j : \R^2 \rightarrow \R^2$ be the
linear function given by $\f_j (r, s) = (r N_j , s M _j)$. We
denote by $R_j $ the number $ R_j := \frac{1}{E_j} N_{j-1} N_j
\bar{B}_j $.  Then we have that:
\[
\mathcal{N} (F_{j+1}, F_j,   \frac{1}{E_j} \underline{\bar{B}}_j,
\underline{F}_j )    = (R_j , R_j)   + \f_j (\mathcal{N}
(F_{j+1}^{(j)} (x_j, y_j))).
\]     \end{remark}

     \subsection{Equisingularity criterions} \label{equi-crit}

Let $F \in \C \{ t, x \}[y]$ be a  Weierstrass polynomial in $y$.
We suppose that  $y$ is an approximate root of $F$, i.e., $F$ is of the form:
\begin{equation}    \label{candidato2}
         F = y^{N} + \sum_{i=2}^N A_{i,{t}} (x) y^{N-i} \in \C \{ x, {t}\}
         [y].
\end{equation}
Set  $F_t (x, y) = F(t, x, y)$ and consider the family of germs $(C_{t}, 0)$ defined
by $F_{t}=0 $.
 We assume that $(x, F_{t})_0  = e_0 > 1$,
for $0 \leq |{t}| \ll 1$.

We give an algorithm to check whether a family of curves $(C_{t},
0)$ of the form (\ref{candidato2}) is equisingular at ${t}=0$ to a
plane branch (irreducible and reduced).  If the answer of the
algorithm is no, then either $(C_0, 0)$ is not analytically
irreducible or $(C_{t}, 0) $ is not equisingular at ${t}=0$. The
proof follows from the discussion in Section \ref{st-line}.
\begin{Alg} \label{cri-1} $\,$

{\bf Step 1}. Set $\mathcal{N}_1$ the Newton polygon of $F = \sum
\a^{(1)}_{r,s} (t) x^r y^s$ (with respect to $(x,y)$).
\begin{enumerate}
\item[(1.a)]  Check that $\mathcal{N}_1$ has only one edge
$\Gamma_1$ with vertices $(e_1 m_1, 0)$ and $(0,e_1 n_1)$ with
$ e_1 \geq 1$, $\mbox{{\rm gcd}} (n_1, m_1) =1$ and $e_1 \geq 1$. If $e_1 =1$
answer yes, otherwise verify that $\min \{n_1, m_1 \} >1$.  Notice that $e_0 = e_1 n_1$.

\item[(1.b)]   Check that  the polynomial $\sum_{k=0}^{e_1}
\a^{(1)}_{e_1 m_1 - k m_1, k n_1} (t) \, z^k  \in \C \{ t \} [z]$
is of the form $(z - \theta_1 (t))^{e_1}$ for some series
$\theta_1 (t) \in \C \{ t\}$ with $\theta_1 (0) \ne 0$.

\end{enumerate}

{\bf Step 2}. If $e_1 > 1$  and conditions {\rm (1.a)}  and {\rm (1.b)} hold,
set $F_2$  for the approximate root of $F$ of degree $n_1$.
Compute the $(y, F_2)$-expansion $F = \sum_{finite}
A^{(2)}_{i_1, i_2} (t, x) \, y^{i_1}  \, F_2^{i_2}$.

\item[(2.a)]  From the data $(n_1, m_1)$ and $e_0$ each triple
\begin{equation} \label{ese}
 (i_0, i_1, i_2) \mbox{ for  } i_0 := \mbox{{\rm ord}}_x
(A^{(2)}_{i_1, i_2} (t, x)), \mbox{ determines } (r,s)  \mbox{
with } r
> 0 \mbox{ and } s < e_1 \mbox{ or } (r, s) = (0, e_1)
\end{equation}
 and the converse also holds.  (This follows by Remark
\ref{llidos} and the method given in the proof of
 Lemma \ref{num}).
For $(r,s)$ and $(i_0, i_1, i_2)$ in  (\ref{ese}) denote $x^{i_0}
y^{i_1} F_2^{i_2}$ by $\mathcal{M}_2(r,s)$.

   \begin{enumerate}
\item[(2.a)] Denote by $\mathcal{N}_2$ for the convex hull of
the set $\cup_{(r,s)}  \{ (r,s)  + \R^2_{\geq 0} \}$,  for $(r,s)$
those of (\ref{ese}). Check that $\mathcal{N}_2$ has only one edge
$\Gamma_2$ with vertices $(e_2 m_2, 0)$ and $(0,e_2 n_2)$ with
$\mbox{{\rm gcd}} (n_2, m_2) =1$ and $e_2 \geq 1$. Verify that $n_2 > 1$ and $m_2 \geq 1$.

\item[(2.b)] If $i_0 =  \mbox{{\rm ord}}_x A^{(2)}_{i_1, i_2} (t,
x)$ we denote by $\a_{r,s} ^{ (2)}(t)$
 the coefficient of $x^{i_0}$ in the expansion of $A^{(2)}_{i_1, i_2} (t, x)$
as a series in $x$, where $(r,s)$ is determined by $(i_0, i_1,
i_2)$ in terms of (\ref{ese}). Set \[
 F^{\Gamma_2}
:= \sum_{k=0}^{e_2} \a^{(2)}_{e_2 m_2 - k m_2, k n_2} (t) \,
\mathcal{M}_2 ( e_2 m_2 - k m_2, k n_2 ) .\] Compute the
approximate root $F_3$ of degree $n_1 n_2$ of $F$. Check that
$F_3^{e_2}$ is of the form   $F_3^{e_2} = F^{\Gamma_2 } + \sum_{
n_2 r + m_2 s  > e_2 n_2 m_2 } \,  \g^{(2)}_{r,s}  (t) \,
\mathcal{M}_2 (r,s) $, for some $ \g^{(2)}_{r,s}  (t) \in \C \{
t\} $. \item[(2.c)] Check that  the polynomial $ \sum_{k=0}^{e_2}
\a^{(2)}_{e_2 m_2 - k m_2, k n_2} (t) \, z^k  \in \C \{ t \} [z] $
is of the form $(z - \theta_2 (t))^{e_2}$,  for some series
$\theta_2 (t) \in \C \{ t\}$ with $\theta_2 (0) \ne 0$.
\end{enumerate}

{\bf Step j > 2}.   If the conditions {\rm (j-1.a), (j-1.b)}  and {\rm (j-1.c)}, corresponding to
{\rm (2.a), (2.b)} and {\rm (2.c)} respectively are verified and $e_{j-1} >1$
compute   $(y, F_2, \dots, F_j)$
expansion of $F$ and check the conditions {\rm (j.a), (j.b)} and  {\rm(j.c)}, corresponding to
{\rm (2.a), (2.b)} and {\rm (2.c)} respectively.

The algorithm stops whenever some condition is not verified,
answering NO, or when all conditions are verified
and $e_g =1$ for some integer $g$,  answering then YES.
\end{Alg}

\begin{remark}
Our criterion extends the one given by A'Campo and Oka in
\cite{AC}. They assume that certain approximate roots of $F_{t}$
do not depend on ${t}$. We do not need this hypothesis. We do not
compute intersection numbers as in Abhyankar's irreducibility
criterion  \cite{Abhyankar} nor resultants as in \cite{GB-G} (see
Subsection \ref{barroso}).
\end{remark}

\begin{example}
We consider $F$ of the form (\ref{candidato2}).
\[
\begin{array}{c}
F:=y^{12} + ( - 6\,x^{3} + 6\,t\,x^{4})\,y^{10} + (15\,x^{6} -
30\,t\,x^{7}
 )\,y^{8} + ( - 20\,x^{9} + (60\,t - 2) \,x^{10})
\,y^{6}
\\
 + (  15\,x^{12}  + ( 6 - 60\,t )\,x^{13}
(-6\,t +  \lambda ) \, x^{14}) \,y^{4}
- x^{16}\,y^{3}
\\
 + ( - 6\,x^{15}  + (- 6 + 30\,t)
\,x^{16}   + (- 2\,\lambda    + 12\,t) \,x^{17})

\,y^{2}
 + x^{19}\,y
+ x^{18} +( 2  - 6\,t) \,x^{19} +  (1
 - 6\,t  + \lambda ) \,x^{20}.
\end{array}
\]
The approximate roots of $F$ of degrees $2$ and $6$ are
$F_2: =
y^2-x^3 + t \, x^4 $ and
$F_3 := F_2^3 -\frac{15}{2}
\,t^2\,x^8\,F_2 +x^{10} (- 1+  20\,t^3\,x^{2} )$.
Notice that
both polynomials depend on the parameter $t$ hence we cannot apply
the equisingularity criterion of \cite{AC}. We check that the
Newton polygon of $F$ has only two vertices $(0, 12)$ and
$(18,0)$. We set $e_1 = 6$, $n_1 =2$ and $m_1 = 3$. The conditions
({1.a}) and ({2.a}) are verified for $F_2$. We compute the
$(y, F_2)$-expansion of $F$ and we find:
\[ \begin{array}{c}
F = F_2^{6} - 15\,t^{2}\,x^{8}\, F_2^{4} + ( - 2 +
40\,t^{3}\,x^{2})\, x^{10} \, F_2^{3} + (\lambda  -
45\,t^{4}\,x^{2}) \, x^{14} \, F_2^{2}  - x^{16}\,y \,F_2 + \\
 ( -
2\,\lambda \,t + 6\,t^{2}  + 24\,t^{5}\,x^{2})\, x^{18}\, F_2 +
t\,x^{ 20} \,y   + ( 1  + (\lambda \,t^{2}  - 4\,t^{3})
\,x^{2}   +  5\,t^{6}\,x^{4})    \, x^{20}
\end{array}
\]
With the notations introduced above we have that $\mathcal{N}_2$
only  two vertices $(0,6)$ and $(4,0)$, hence $e_2 = 2$, $n_2 = 3$
and $m_2 = 2$. We have that
$
F^{\Gamma_2}  = \mathcal{M}_2 (0,6) - 2 \mathcal{M}_2 (3,2) +
\mathcal{M}_2 (0,4)$,
where $ \mathcal{M}_2 (0,6) = F_2^6 $,  $\mathcal{M}_2 (3,2) =
x^{10} \, F_2^3 $ and $ \mathcal{M}_2 (0,4)= x_1^{20}$. We check
that the conditions ({2.b}) and ({2.c}) are satisfied. We
compute the $(y, F_2, F_3)$ expansion of $F$, for $F_3$ the
approximate root of degree $6$ of $F$. We obtain that
\[ \begin{array}{c}
F= F_3^2
 + (   \lambda -  \frac {405}{4} \,t^{4}\,x^{2}) \, x^{14}
 \,F_2^{2}  - x^{16}\,y  \, F_2
 +
 \\
 ( -2\,\lambda \,t\, - 9\,t^{2} + 324\,t^{5}\,x^{2} ) \, x^{18} \,F_2
 + t \,x^{20}\,y + ( \lambda
\,t^{2} + 36\,t^{3}  - 405\,t^{6}\,x^{2} ) \,x^{22}.
\end{array}
\]
In order to compute the polygon $\mathcal{N}_3$ we consider the
leading terms in the expansion above and we use the method of Lemma \ref{num}.
We have that $x^{i_0} y^{i_1} F_2^{i_2} F_3^{i_3} =
\mathcal{M}_2 (r_1, s_1) F_3^{i_3} $  where
$s_1 = i_2 $ and $r_1 = i_0 n_1 +  m_1 ( -e_0 + i_1 + n_1 i_2 + n_1 n_2 i_3)$.
  We have then that  $\mathcal{M}_2 (r_1, s_1) F_3^{i_2} = \mathcal{M}_3 (r, s)$ where $s = i_2$ and
$r = n_2 r_1 + m_1 ( -e_1  + s_1 + n_2 i_3)$. For instance,
 we have that  $x^{14}\,
F_2^2 = \mathcal{M}_2 (4,2) = \mathcal{M}_3 (4,0)$ and
$x^{16} y F_2 = \mathcal{M}_2 ( 5,1) = \mathcal{M}_3 (5,0)$.
We distinguish two cases in terms of the constant $\lambda \in
\C$.
\begin{enumerate}
\item[(a)]
 If $\lambda \ne 0$ then we check that $\mathcal{N}_3$
is a polygon with vertices $(0,2)$ and $(4,0)$. We have that $e_3
= 2$, $n_3 =1$ and $m_3= 2$ hence  $\{F_t =0 \}$ is not
equisingular at $t = 0$. \item[(b)] If $\lambda =0$ then
$\mathcal{N}_3$ is a polygon with vertices $(0,2)$ and $(5,0)$. We
have that $e_3 =1$ and the conditions (2.a) (2.b) and (2.c) are verified,
hence $\{F_t =0 \}$ is equisingular at $t=0$.
\end{enumerate}
\end{example}

\subsubsection{Equisingularity criterion by jacobian Newton
polygons}     \label{barroso}

Let $f \in \C\{ x\} [y]$ be a Weierstrass polynomial. The {\em
jacobian Newton polygon} of $f$ with respect to the line $ \{ x =0 \}$ is
the Newton polygon of $\mathcal{J}_f(s, x):= \mbox{\rm Res}_y (s -
f, \frac{\partial{f}}{\partial y}) \in \C \{ x, s \}$, where $\mbox{\rm Res}_y$
denotes the resultant with respect to $y$. The jacobian Newton
polygon appears in more general contexts related to invariants of
equisingularity (see \cite{T1}). Garc\'\i a Barroso
and Gwo\'zdziewicz have proved that if $f' \in \C \{ x \} [y] $ is
irreducible and
               $\mathcal{J}_f(s, x) =  \mathcal{J}_{f'}(s, x)$ then $f$ is irreducible.
They have given two methods which characterize jacobian polygons
of plane branches among other Newton polygons by a finite number
of combinatorial operations on the polygons (see \cite{GB-G},
Theorem 1, 2 and 3). The following algorithm is consequence of
their work.
\begin{Alg} \label{cri-2}
Imput: A family $F_{t}(x, y)$ of the form (\ref{candidato2}).

\begin{enumerate}
\item[(a)]
 Compute $\mathcal{J}_{F_{{t}}} (s, x)$.

\item[(b)] Compute the
Newton polygon $\mathcal{N}_{t}$  of     $\mathcal{J}_{F_{{t}}}
(s, x)$.  Check that $  \mathcal{N}_{t}= \mathcal{N}_0$.

\item[(b)] Test if    $\mathcal{N}_0$  is a jacobian Newton polygon of
a plane branch by using Theorem 2 or 3 in \cite{GB-G}.
\end{enumerate}
If all the steps of the algorithm give a positive answer then
$F_{t}=0$ is equisingular at ${t}=0$ to a plane branch.
\end{Alg}

\section{Multi-semi-quasi-homogeneous deformations} \label{msqh}

In this Section we introduce a class of (non equisingular)
deformations of a plane branch $(C,0)$ and we study some of its basic properties
which are essential for the applications in the real case (see \cite{GP-R}).

We keep notations of Section \ref{toric-res}. The resolution is
described in terms of a fixed sequence $f_1, \dots,f_g$ of
semi-roots. For simplicity we assume that  $\theta_1, \dots ,
\theta_g =1$ (see Remark \ref{normalized}). We introduce the
following notations:
\begin{notation} \label{notation-res}
For $j=1, \dots,g$ we set:
\begin{enumerate}
\item[(i)]  $ \Gamma_j = [(m_{j} e_{j}, 0) , ( 0, n_j e_j)]$ the
compact edge of the local Newton polygon of $f^{(j)} (x_j, y_j)$
(see (\ref{fg})).

\item[(ii)] $\D_j$ the triangle bounded by the Newton polygon of
$f^{(j)} (x_j, y_j)$ and the coordinate axis; we denote by $\D_j^-
$ the set $\D_j^- = \D_j \setminus \Gamma_j$.

\item[(iii)]
 Let $\w_j: \D_{j} \cap \Z^2 \rightarrow \Z$ be defined by
$
 \w_j
(r, s) = e_{j}  ( e_{j} n_{j} m_{j} - r n_{j} - s m_{j} ).
$
\end{enumerate}
\end{notation}
The symbol $\tr_j$ denotes the parameters $(t_j, \dots, t_g)$ for
any $1 \leq j \leq g$. We consider sequences of multiparametric
deformations $C_{\tr_1}, \dots, C_{\tr_g}$  of $(C,0)$ defined by
$P_{\tr_1}, \dots, P_{\tr_g}$ of the form (see Definition
\ref{opq}):
\begin{equation} \label{pg}
\left\{
\begin{array}{cccccl}
P_{\tr_g} & := &  F & +  & \displaystyle{\sum_{(r,s) \in \D_{g}^-
\cap \Z^2}} \, a_{r,s}^{(g)} (t_g)
 & M_{g} (r,s)
\\
\dots & \dots & \dots & \dots & \dots & \dots
\\
 P_{\tr_1 } & := &  P_{\tr_2} & +  & \displaystyle{\sum_{(r,s) \in \D_{1}^- \cap
\Z^2}} \,  a_{r,s}^{(1)} (t_1) & M_1(r,s),
\end{array}
\right.
\end{equation}
where  $a_{r,s}^{(j)} (t_j) \in \C \{ t_j \} $  for $(r,s) \in
\D_j^- \cap \Z^2$ and $j=1, \dots, g$.
 Notice that $P_{\tr_1}$
determines any of the terms $P_{\tr_j}$ for $1 < j \leq g$, by
substituting
 $t_{1} = \dots = t_{j-1}
=0$ in $P_{\tr_1}$.
 The multiparametric deformation $C_{\tr_1}$ is
 {\em multi-semi-quasi-homogeneous}  (msqh) if in addition
$a^{(j)}_{r,s} =  A^{(j)}_{r,s} \, t_j^{\w_j (r,s)}$, for $1\leq j
\leq g$ and $(r,s) \in \D_j^- \cap \Z^2$, where $\w_j(r,s) \in
\Z_{\geq 0}$ is defined in Notation \ref{notation-res},
$A^{(j)}_{r,s} \in \C$ and $ A^{(j)}_{0,0} \ne 0$, for $j=1,
\dots, g$.

\begin{remark}
In the real case we have studied the topological types of the
msqh-deformations with real part with the maximal number of
connected components. The
hypothesis of being msqh is related in that paper to the study of
the assymtotic scales of the ovals when the parameters tend to
zero (joint work with Risler \cite{GP-R}).
\end{remark}

We denote by $C_{l, \tr}^{(j)} \subset Z_j$ the strict transform
of $C_{l, \tr}$ by the composition of toric maps $\p_{1} \circ
\cdots \circ \p_{j-1}$ and by $P_{\tr_l}^{(j)} (x_j, y_j)$
the polynomial defining
$C_{\tr_l}^{(j)}$ in the coordinates $(x_j, y_j)$
for $2 \leq j \leq l \leq g$. These notations are
analogous to those used for $C$ in Section \ref{toric-res}, see
(\ref{strict-j}).

\begin{proposition} \label{perturbation}
If $1 \leq j < l \leq g$ the curves $C_{\tr_l}^{(j)}$ and
$C^{(j)}$ meet the exceptional divisor of $\p_{1} \circ \cdots
\circ \p_{j-1}$ only at the point $o_j \in \{ x_j =0 \} $ and with
the same intersection multiplicity $e_{j-1}$.
\end{proposition}
{\em Proof.} If $j=1$ we have that $f$ and $P_{\tr_l}$ have the
same Newton polygon and moreover the symbolic restrictions of
these two polynomials to the compact face $\Gamma_1$ of the Newton
polygon coincide by Lemma \ref{newton}. If $j >1$ we show the
result by induction using Remark \ref{newton-rem}. $\Box$

\begin{proposition} \label{perturbation2}
If $ 1 < j \leq g$
then $\{ x_j = 0 \}$ is the only irreducible component of the
exceptional divisor of $\p_{1} \circ \cdots \circ \p_{j-1}$ which
intersects $C_{\tr_j}^{(j)}$ at $e_{j-1}$ points counted with
multiplicity. More precisely, we have that:
\begin{enumerate}
\item[(i)]The symbolic restriction of $P^{(j)}_{\tr_{j+1}} (x_{j},
y_{j}) $ to the edge $\Gamma_j$ of its local Newton polygon is of
the form: $ \a_j \prod_{s=1}^{e_{j}} \left( y_{j} ^{n_j} - (1 +
\g_s^{(j)} t_{j+1}^{e_{j+1} m_{j+1} } ) x_{j}^{m_j} \right)$, with
$ \a_j, \g_s^{(j)} \in \C^*$, for $s=1, \dots, e_{j}$.
 \item[(ii)] The points of intersection of  $\{ x_{j+1} = 0 \} $ with
$C_{\tr_{j+1}}^{(j+1)}$ are those with coordinates $x_{j+1} =0$
and
\begin{equation} \label{intersection2}
u_{j+1} = (1 + \g_s^{(j)} t_{j+1}^{e_{j+1} m_{j+1}})^{-1} , \mbox{
for } s=1, \dots, e_{j}.
\end{equation}
\end{enumerate}
\end{proposition}
{\em Proof.} If $j=2$ we have that the terms of the expansion of
$P_{\tr_2}$ which have exponents on the compact face of the Newton
polygon  $\mathcal{N} (f) $ are $f$ and ${\mathcal{M}}_2 (0, s)$
for $0 \leq s < e_1$ by Lemma \ref{newton}. By Proposition
\ref{comb2} we have that ${\mathcal{M}}_2 (0,s) = x_1^{m_1 (e_1
-s)} f_2^s$. The restriction of the polynomial $ f +
\sum_{s=0}^{e_1 - 1} A_{(0,s)}^{(2)} t_2^{e_2 m_2 (e_1 -s) }
x_1^{m_1 (e_1 -s)} f_2 ^s$ to the face $\Gamma_j$ is equal to:
\begin{equation} \label{expression-2}
(y_1^{n_1} - x_1^{m_1})^{e_1} + \sum_{s=0}^{e_1 - 1}
A_{(0,s)}^{(2)} (t_2^{e_2 m_2} x_1^{m_1})^{e_1 -s} (y_1 -
x_1^{m_1})^s .
\end{equation}
Let us consider the polynomial $Q_1 (V_1, V_2) := V_1^{e_1} +
\sum_{s=0}^{e_1 - 1} A_{(0,s)}^{(2)} V_1^{s} V_2^{e_1 - s}$. By
hypothesis $A_{(0,0)}^{(2)} \ne 0$ hence the homogeneous
polynomial $Q_1$ factors as: $ Q_1 (V_1, V_2) = \prod_{s=1}^{e_1}
(V_1 - \g_s^{(1)} V_2) $ for some $ \g_s^{(1)}\in \C^*$. The
expression (\ref{expression-2}) is of the form: $ Q_1 (y_1^{n_1} -
x_1^{m_1}, t_2^{m_2 e_2} x_1^{m_1}) = \prod_{s=1}^{e_1} (y^{n_1} -
(1 + \g_s^{(1)}  t_2^{m_2 e_2}) x_1^{m_1})$. This proves the first
assertion in this case. The second follows  from this by the
discussion of Section \ref{toric-res}.

If $j >2$ we deduce by induction, by using Remarks
\ref{newton-rem} and \ref{units}, that the restriction of
$P_{\tr_{j+1}}^{(j-1)}$ to the compact face of $\mathcal{N}
(f^{(j-1)} )$ is of the form $(y_{j}^{n_j} - x_1^{m_j})^{e_{j}} +
\sum_{s=0}^{e_{j} - 1} A_{(0,s)}^{(j+1)} (t_{j+1}^{e_{j+1}
m_{j+1}} x_j^{m_{j+1}})^{e_{j} -s} (y_j ^{n_j} - x_j^{m_j})^s$.
The result  follows by the same argument. $\Box$.

\subsection{Milnor number and generic msqh-smoothings}

If $(D, 0) \subset (\C^2, 0)$ is the germ of a plane curve
singularity, defined by $h =0$, for $h \in \C \{ x, y \} $
reduced, we denote by $\mu(h)_0$ or by $\mu (D)_0$ the {\em Milnor
number} $ \dim_\C \C \{ x, y \} / (h_x, h_y) $. We have the
following formula (see \cite{Risler2} and \cite{Z}):
\begin{equation}\label{m-milnor}
 \mu(h)_0 = (h, \frac{\partial h}{\partial y})_0 - (h, x)_0 +
1.
\end{equation}
The Milnor number  of the plane branch $(C, 0)$ expresses also in terms of the generators of the
semigroup $\Lambda_C$ with respect to the coordinate line $\{ x=0
\}$ (see \cite{Merle,Barroso} and \cite{Z}) by using:
\begin{equation} \label{mil-semi}
(f, \frac{\partial f}{\partial y})_0   = \sum_{j=1}^g (n_j -1)
{\bar{b}}_j.
\end{equation}

    \begin{definition}
We say that the deformation $P_{\tr_1}$ of the plane branch $(C,
0)$ is {\em generic} if the numbers $\{ \g_s^{(j)} \}_{s =1}^{e_j}
$ appearing in Proposition \ref{perturbation2} are all distinct,
for $1 < j \leq g$.
    \end{definition}

    The following Proposition provides a geometrical
    incarnation in terms of the
    sequence of generic msqh-deformations
    of the Milnor's formula $\mu(C)_0 =
    \frac{1}{2} \delta(C)_0$, for $\delta(C)_0$
    the delta invariant of  $(C,0)$ (see \cite{Wall} , \cite{Casas} Ex. 5.6, see also \cite{Janus}).
\begin{proposition} \label{milnor}
Let $P_{\tr_1}$ be a generic msqh-deformation of a plane branch
$(C, 0)$,  then we have that \[  \mu (C)_{0} = \displaystyle{
\sum_{ j=1}^g ( \mu ( C_{\tr_{j+1}}^{(j)} )_{o_{j}} + e_j -1 )}
\]
\end{proposition}
{\em Proof.} We prove the result by induction on $g$. If $g =1$
the assertion is trivial. We suppose the assertion true for
branches with $g-1$ characteristic exponents with respect to some
system of coordinates. By Proposition \ref{mod} and the induction
hypothesis it is easy to check that $ \mu (C^{(2)} )_{o_2} = {
\sum_{ j=2}^g ( \mu ( C_{\tr_{j+1}}^{(j)} )_{o_{j}} + e_j -1 )}$.

 By Proposition \ref{perturbation2} and the
definition of generic msqh-deformation, we have that the curve
$C_{\tr_2}^{(1)}$, defined by the polynomial $P_{\tr_2}^{(1)}
(x_1, y_1)$, is non-degenerate with respect to its Newton
polygon. By (\ref{m-milnor})  we have that $
   \mu ( C_{\tr_2}^{(1)} )_{o_1} = e_0 b_1 - e_0 - b_1 +1$. By (\ref{mil-semi}) we have that: $
 \mu (C)_{o_1}  -  \mu (C^{(2)} )_{o_2}  =  (f,\frac{\partial f}{\partial y})_0    - (f^{(2)}, \frac{\partial f ^{(2)}}{\partial
y_2} )_{o_2}                      + e_1 - e_0 $.   The assertion
holds if and only if $ (f,\frac{\partial f}{\partial y})_0    -
(f^{(2)}, \frac{\partial f ^{(2)}}{\partial y_2} )_{o_2} = b_1
(e_0 -1)$. Using (\ref{mil-semi}) and Lemma \ref{j-semi} we verify
that $ (f,\frac{\partial f}{\partial y})_0    - (f^{(2)},
\frac{\partial f ^{(2)}}{\partial y_2} )_{o_2}$ is equal to:
${\sum_{j=1}^g} (n_k -1) {\bar{b}}_j - {\sum_{j=2}^g} (n_k -1)
\left( {\bar{b}}_j - {\bar{b}}_1 e_0 /e_{j-1} \right)   =  (n_1 -1
) {\bar{b}}_1 + {\sum_{j=2}^g} {\bar{b}}_1 e_0/e_{j-1} =
 \left(   n_1 -1 + (n_2 -1) n_1 + \cdots +
(n_g -1) n_1 \cdots n_{g-1} \right) {\bar{b}}_1 =  (e_0-1) b_1$.
$\Box$

\begin{corollary}
\[ \mu (C)_0 = 2 (\sum_{j=0}^{g-1}  \left( \# (
\stackrel{\circ}{\D}_j \cap \Z^2 ) + e_{j+1} -1 \right) ) . \]
\end{corollary}

{\bf Acknowledgement.} The author is grateful to J.-J. Risler, B. Teissier and  
P. Popescu-Pampu for useful discussions on
this topic and to Arkadiusz P\l oski and Evelia Garc\'\i a Barroso 
for pointing out a problem in the proof of Theorem 3.8 in 
a previous version of this paper. 

 {\small


\begin{thebibliography}{GGGG}


\bibitem[A'C-Ok]{AC}{\sc A'Campo, N., Oka, M.}, Geometry of plane curves via Tschirnhausen resolution tower,
{\em Osaka J. Math.}, {\bf 33}, (1996), 1003-1033.



\bibitem[Abh2]{A-expansions} {\sc Abhyankar, S.S.}, {\em Lectures on expansion techniques in algebraic geometry.
Notes by Balwant Singh.} Tata Institute of Fundamental Research Lectures on Mathematics and Physics, 57.
Tata Institute of Fundamental Research, Bombay, 1977.


\bibitem[Abh3]{A-Kyoto} {\sc Abhyankar, S.S.,} On the semigroup of a meromorphic curve. I.
{\em Proceedings of the International Symposium on Algebraic Geometry (Kyoto Univ., Kyoto, 1977)},  pp. 249--414, Kinokuniya Book Store, Tokyo, 1978.


\bibitem[Abh4]{Abhyankar} {\sc Abhyankar, S.S.},
Irreducibility criterion for germs of analytic functions of two complex variables.
{\em Adv. Math.} {\bf 74} (1989), no. 2, 190--257.



\bibitem[A-M]{Moh}{\sc Abhyankar, S.S.,  Moh, T.}, Newton-Puiseux
  Expansion
and Generalized Tschirnhausen Transformation I-II, {\em  J.Reine
Angew. Math.},  {\bf 260}. (1973), 47-83;  {\bf
  261}. (1973), 29-54.


\bibitem[A-M2]{Moh2}  {\sc Abhyankar, S.S.,  Moh, T.},  Embeddings of the Line in the Plane. {\em
J.Reine Angew. Math.}, {\bf 276} (1975), 148-166.


\bibitem[As-B]{Assi} {\sc Assi, A., Barile M.}, Effective construction of irreducible curve singularities.
{\em Int. J. Math. Comput. Sci. } {\bf 1}  (2006),  no. 1, 125--149.

\bibitem[C]{Campillo}{\sc  Campillo, A., }
 {\em Algebroid Curves in positive characteristic},  Lecture Notes
 in Mathematics (813),
  Springer,  1980.



\bibitem[Ca]{Casas} {\sc Casas-Alvero, E.},
{\em Singularities of plane curves}
London Mathematical Society Lecture Note Series 276,
Cambridge University Press,
Cambridge,
2000.



   \bibitem[C-M1]{Cossart-Saskatoon}   {\sc  Cossart, V., Moreno-Soc\' \i as, G. }
Irreducibility criterion: a geometric point of view.
{\em  Valuation theory and its applications, Vol. II (Saskatoon, SK, 1999)},  27--42,
Fields Inst. Commun., 33, Amer. Math. Soc., Providence, RI, 2003.



\bibitem[C-M2]{Cossart}{\sc  Cossart, V., Moreno-Soc\' \i as, G. }
Racines approch\'ees, suites g\'en\'eratrices, suffisance des jets, {\em
Annales de la Facult\'e des Sciences de Toulouse} Vol XIV, no 3,
2005, pp. 353-394.

\bibitem[GB]{Barroso} {\sc Garc\'\i a Barroso, E.R.}, Sur les courbes polaires d'une courbe plane r\'eduite.
{\em Proc. London Math. Soc. } {\bf 3} 81 (2000),  no. 1, 1-28.

\bibitem[GB-G]{GB-G} {\sc Garc\'\i a Barroso, E.R.,   Gwo\'zdziewicz},  Characterization of jacobian Newton polygons of plane branches and new
criteria of irreductibility.
ArXiv: 0805.4257[math.AG], to appear in    {\em
  Ann. Inst. Fourier (Grenoble)}.



\bibitem[G-T]{Rebeca}{\sc  Goldin, R., Teissier, B.,}
Resolving singularities of plane analytic branches with one toric
morphism, {\em Resolution of Singularities, A research textbook in
tribute to Oscar Zariski}. Edited by H. Hauser, J. Lipman, F.Oort
and A. Quiros. Progress in Mathematics No. 181,
Birkh\"auser-Verlag, 2000, 315-340.

\bibitem[GP]{GP} {\sc Gonz\'alez P\'erez, P.D.},  Toric embedded
  resolutions of quasi-ordinary hypersurface singularities.  {\em
  Ann. Inst. Fourier (Grenoble)} {\bf 53} (2003),  no. 6, 1819--1881.


\bibitem[GP-R]{GP-R} {\sc Gonz\'alez P\'erez, P.D.; Risler, J.-J.},
Multi-Harnack smoothings of real plane branches, arXiv:0808.0157v1 [math.AG].



\bibitem[G-P]{Ploski} {\sc Gwo\' zdziewicz, J., Ploski, A.},  On the
Approximate Roots of Polynomials. {\em Annales Polonici
Mathematici}, {\bf LX} 3,  (1995), 199-210.




\bibitem[G]{Janus} {\sc Gwo\' zdziewicz, J.}, Kouchnirenko type formulas for local invariants of
             plane analytic curves, arXiv:0707.3404v1 [math.AG]

\bibitem[L-Ok]{LO} {\sc
L\^e D.T. and  Oka, M.;} On resolution complexity of plane curves.
{\em Kodai Math. J.} {\bf 18}  (1995),  no. 1, 1--36.






\bibitem[M]{Merle}
{\sc Merle, M.,}      Invariants polaires des courbes planes.
{\em Invent. Math.} 41 (1977), no. 2, 103--111.




\bibitem[Ok1]{Oka1}{\sc Oka, M.}, Geometry of plane curves via toroidal resolution,
{\em Algebraic Geometry and Singularities}, Edited by A. Campillo
L\'opez and L. Narv\'aez Macarro. Progress in Mathematics No. 139,
Birkh\"auser, Basel, 1996.

\bibitem[Ok2]{Oka2}{\sc Oka, M.}{\em  Non-degenerate complete intersection singularity},
Actualit\'es Math\'ematiques, Hermann, Paris, 1997.


\bibitem[Pi]{Pinkham} {\sc Pinkham, H.}, Courbes planes ayant une seule place a l'infini. {\em S\'eminaire sur les Singularit\'es
des surfaces.} Centre de Math\'ematiques de l'\'Ecole Polytechnique. Ann\'ee 1977-1978.

\bibitem[PP]{PP} {\sc  Popescu-Pampu, P.}, Approximate roots. {\em Valuation theory and its applications, Vol. II (Saskatoon, SK, 1999)},  Fields Inst. Commun., 33, Amer. Math. Soc., Providence, RI, 2003, 285-321.


\bibitem[R]{Risler2}
{\sc Risler, J-J.,}           Sur l'id\'eal jacobien d'une courbe plane,
{\em Bulletin de la S. M. F.}  {\bf 99} (1971), 305--311.


\bibitem[T1]{T1}{\sc  Teissier, B.,}
 {\sc  Teissier, B.,}  Vari\'et\'es polaires. I. Invariants polaires des singularit\'es d'hypersurfaces.
{\em  Invent. Math.} {\bf 40}  (1977), no. 3, 267--292.


\bibitem[T2]{T}{\sc  Teissier, B.,} The monomial curve and its
  deformations, Appendix in \cite{Z}

    \bibitem[T3]{Tcurvas}{\sc  Teissier, B.,}
Complex curve singularities: a biased introduction.
{\em Singularities in geometry and topology},  825--887, World Sci. Publ., Hackensack, NJ, 2007.


\bibitem[W]{Wall}
 {\sc Wall, C. T. C.} {\em Singular points of plane curves.}
 London Mathematical Society Student Texts,  63. Cambridge University Press, Cambridge, 2004.


\bibitem[Z1]{Z-inversion}{\sc  Zariski, O.,} Studies in
  equisingularity III, Saturation of local rings and equisingularity
on Algebraic Varieties, {\em Amer. J. Math.},  {\bf  90},  1968 961--1023.

\bibitem[Z2]{Z}{\sc  Zariski, O.,}
{\em Le probl\`eme des modules pour les branches planes},
Hermann, Paris, 1986.

\end{thebibliography}
\end{document}